\documentclass{amsart}
\usepackage{amssymb}
\usepackage{amscd}
\usepackage{verbatim}
\usepackage{epsfig}
\begin{document}
\newcommand\Mwhere{\ \text{where}\ }
\newcommand\Mand{\ \text{and}\ }
\newcommand\Mwith{\ \text{with}\ }
\newcommand\Mfor{\ \text{for}\ }
\newcommand\Mst{\ \text{such that}\ }
\newcommand\Mor{\ \text{or}\ }
\newcommand\Mif{\ \text{if}\ }
\newcommand\Miff{\ \text{iff}\ }
\newcommand\Mthen{\ \text{then}\ }
\newcommand\nin{\notin}
\newcommand\identity{\operatorname{id}}
\newcommand\Id{\operatorname{Id}}
\newcommand\Real{\mathbb{R}}
\newcommand\pos{\Real^+}
\newcommand\Rnp{\Real\setminus\{0\}}
\newcommand\nzero{\setminus\{0\}}
\newcommand\Cx{\mathbb{C}}
\newcommand\Cxp{\Cx^+}
\newcommand\Cxm{\Cx^-}
\newcommand\Nat{\mathbb{N}}
\newcommand\halfNat{{\frac{1}{2}}\mathbb{N}}
\newcommand\intgr{\mathbb{Z}}
\newcommand\im{\operatorname{Im}}
\newcommand\re{\operatorname{Re}}
\newcommand\sign{\operatorname{sign}}
\newcommand\codim{\operatorname{codim}}
\newcommand\End{\operatorname{End}}
\newcommand\Ker{\operatorname{Ker}}
\newcommand\Hom{\operatorname{Hom}}
\newcommand\tr{\operatorname{tr}}
\newcommand\Tr{\operatorname{Tr}}
\newcommand\ideal{{\mathcal I}}
\newcommand\Span{\operatorname{span}}
\newcommand\image{\operatorname{image}}
\newcommand\Range{\operatorname{Ran}}
\newcommand\Graph{\operatorname{graph}}
\newcommand\slim{\operatornamewithlimits{s-lim}}
\newcommand\sll{\operatorname{sl}}
\newcommand\sol{\operatorname{so}}
\newcommand\SL{\operatorname{SL}}
\newcommand\SO{\operatorname{SO}}
\newcommand\On{\operatorname{O}}
\newcommand\pa{\partial}
\newcommand\Rn{\Real^n}
\newcommand\Rm{\Real^m}
\newcommand\RN{\Real^N}
\newcommand\RtN{\Real^{2N}}
\newcommand\RM{\Real^M}
\newcommand\sphere{\mathbb{S}}
\newcommand\Sn{\sphere^{n-1}}
\newcommand\Sm{\sphere^{m-1}}
\newcommand\Snp{\sphere^n_+}
\newcommand\Smp{\sphere^m_+}
\newcommand\SN{\sphere^{N-1}}
\newcommand\SNp{\sphere^N_+}
\newcommand\circlep{\sphere^1_+}
\newcommand\Phom{P_{h}}
\newcommand\Shom{S_{h}}
\newcommand\distance{\operatorname{dist}}
\newcommand\cl{\operatorname{cl}}
\newcommand\interior{\operatorname{int}}
\newcommand\Fa{\operatorname{Fa}}
\newcommand\ff{\operatorname{ff}}
\newcommand\mf{\operatorname{mf}}
\newcommand\cf{\operatorname{cf}}
\newcommand\scf{\operatorname{sf}}
\newcommand\lf{\operatorname{lf}}
\newcommand\rf{\operatorname{rf}}
\newcommand\indfam{{\mathcal K}}
\newcommand\fraka{{\mathfrak a}}
\newcommand\calA{{\mathcal A}}
\newcommand\calB{{\mathcal B}}
\newcommand\calR{{\mathcal R}}
\newcommand\calO{{\mathcal O}}
\newcommand\calJ{{\mathcal J}}
\newcommand\calM{{\mathcal M}}
\newcommand\calN{{\mathcal N}}
\newcommand\calX{{\mathcal X}}
\newcommand\calF{{\mathcal F}}
\newcommand\calG{{\mathcal G}}
\newcommand\calT{{\mathcal T}}
\newcommand\calC{{\mathcal C}}
\newcommand\calCt{{\tilde {\mathcal C}}}
\newcommand\calCL{{\mathcal C}_{\text L}}
\newcommand\calCR{{\mathcal C}_{\text R}}
\newcommand\Cinf{{\mathcal C}^{\infty}}
\newcommand\dist{{\mathcal C}^{-\infty}}
\newcommand\dCinf{\dot\Cinf}
\newcommand\ddist{\dot\dist}
\newcommand\Cj{{\mathcal C}^j}
\newcommand\Linf{L^{\infty}}
\newcommand\phg{{\text{phg}}}
\newcommand\bcon{{\mathcal A}}
\newcommand\bconc{{\mathcal A}_{\text{phg}}}
\newcommand\Sch{{\mathcal S}}
\newcommand\temp{\Sch^{\prime}}
\newcommand\Diff{\operatorname{Diff}}
\newcommand\Diffb{\operatorname{Diff}_{\text{b}}}
\newcommand\Diffc{\operatorname{Diff}_{\text{c}}}
\newcommand\Diffsc{\operatorname{Diff}_{\text{sc}}}
\newcommand\DiffI{\operatorname{Diff}_{\text{I}}}
\newcommand\DiffIq{\operatorname{Diff}_{\text{I},q}}
\newcommand\sing{\text{sing}}
\newcommand\reg{\text{reg}}
\newcommand\supp{\operatorname{supp}}
\newcommand\ssupp{\operatorname{sing\ supp}}
\newcommand\csupp{\operatorname{cone\ supp}}
\newcommand\esupp{\operatorname{ess\ supp}}
\newcommand\Fr{{\mathcal F}}
\newcommand\Frinv{\Fr^{-1}}
\newcommand\bop{{\mathcal B}}
\newcommand\spec{\operatorname{spec}}
\newcommand\pspec{\spec_{pp}}
\newcommand\cspec{\spec_{c}}
\newcommand\FIO{{\mathcal I}}
\newcommand\SP{\operatorname{RC}}
\newcommand\RC{\operatorname{RC}}
\newcommand\Symc{S_c}
\newcommand\Symca{S_c^{\alpha}}
\newcommand\Symczero{S_c^{0,...,0}}
\newcommand\sci{{}^{\text{sc}}}
\newcommand\sct{\sci T^*}
\newcommand\scdt{\sci \dot T^*}
\newcommand\dS{\dot S^*}
\newcommand\dT{\dot T^*}
\newcommand\dSreg{\dot\Sigma_{\text reg}}
\newcommand\scct{\sci\bar{T}^*}
\newcommand\Csc{C_{\text{sc}}}
\newcommand\SNpscd{(\SNp)^2_{\text{sc}}}
\newcommand\scdiag{\Delta_{\text{sc}}}
\newcommand\projscl{\pi^L_{\text{sc}}}
\newcommand\projscr{\pi^R_{\text{sc}}}
\newcommand\scHL{\sci H^{2,0}_{|\zeta|^2-\lambda^2}}
\newcommand\scHrg{\sci H^{2,0}_{\sqrt{g}}}
\newcommand\Hsc{H_{\text{sc}}}
\newcommand\WF{\operatorname{WF}}
\newcommand\WFp{\operatorname{WF^{\prime}}}
\newcommand\WFsc{\operatorname{WF}_{\text{sc}}}
\newcommand\WFscp{\operatorname{WF_{sc}^{\prime}}}
\newcommand\WFC{\operatorname{WF}_C}
\newcommand\WFCi{\operatorname{WF}_{C_i}}
\newcommand\elliptic{\operatorname{ell}}
\newcommand\Psop{\operatorname{\Psi}}
\newcommand\Psiscrs{\operatorname{\Psi_{sc}^{-2,\infty}}}
\newcommand\Psiscr{\operatorname{\Psi_{sc}^{-2,0}}}
\newcommand\Psiscrm{\operatorname{\Psi_{sc}^{0,2}}}
\newcommand\PsiscHam{\operatorname{\Psi_{sc}^{2,0}}}
\newcommand\Psisci{\operatorname{\Psi_{sc}^{*,*}}}
\newcommand\Psiscid{\operatorname{\Psi_{sc}^{0,0}}}
\newcommand\Psiscis{\operatorname{\Psi_{sc}^{0,\infty}}}
\newcommand\Psiscsi{\operatorname{\Psi_{sc}^{-\infty,0}}}
\newcommand\Psiscs{\operatorname{\Psi_{sc}^{-\infty,\infty}}}
\newcommand\Psiscalg{\operatorname{\Psi_{sc}^{\infty,-\infty}}}
\newcommand\nullHam{{\mathcal N}}
\newcommand\charD{\Sigma_{\Delta-\lambda^2}}
\newcommand\charLap{\Sigma_{\Delta-\lambda}}
\newcommand\Snl{\Sn_{\lambda}}
\newcommand\SNl{\SN_{\lambda}}
\newcommand\gammat{\tilde\gamma}
\newcommand\gammasc{\gamma}
\newcommand\Tau{\mathcal{T}}
\newcommand\taut{\tilde\tau}
\newcommand\taub{\bar\tau}
\newcommand\Nout{N^+_{\lambda}}
\newcommand\Nin{N^-_{\lambda}}
\newcommand\Nio{N^{\pm}_{\lambda}}
\newcommand\El{E_{\lambda}}
\newcommand\Elt{\tilde E_{\lambda}}
\newcommand\Eil{E^i_{\lambda}}
\newcommand\Ejl{E^j_{\lambda}}
\newcommand\Eajl{E^{\alpha_j}_{\lambda}}
\newcommand\Eilt{\tilde E^i_{\lambda}}
\newcommand\Np{N^+}
\newcommand\Nm{N^-}
\newcommand\Npm{N^{\pm}}
\newcommand\Fin{F^-(\lambda)}
\newcommand\Fini{F^-_i(\lambda)}
\newcommand\Fout{F^+(\lambda)}
\newcommand\Fouti{F^+_i(\lambda)}
\newcommand\Foutj{F^+_j(\lambda)}
\newcommand\Rout{R^+_{\lambda}}
\newcommand\Routl{R^+_{\lambda^2}}
\newcommand\Routsgnl{R^{\sign\lambda}_{\lambda^2}}
\newcommand\Rin{R^-_{\lambda}}
\newcommand\Rinl{R^-_{\lambda^2}}
\newcommand\Rinsgnl{R^{-\sign\lambda}_{\lambda^2}}
\newcommand\Rio{R^{\pm}_{\lambda}}
\newcommand\Riol{R^{\pm}_{\lambda^2}}
\newcommand\Roi{R^{\mp}_{\lambda}}
\newcommand\Roil{R^{\mp}_{\lambda^2}}
\newcommand\Riob{R^{\pm}}
\newcommand\Roib{R^{\mp}}
\newcommand\Tio{T^{\pm}}
\newcommand\Tiob{T^{\pm}_{\ff}}
\newcommand\Toi{T^{\mp}}
\newcommand\Toib{T^{\mp}_{\ff}}
\newcommand\TIiob{T_I^{\pm}}
\newcommand\Rinb{R^-}
\newcommand\Rinbsgnl{R^{-\sign\lambda}}
\newcommand\Tin{T^-}
\newcommand\Tinb{T^-_{\ff}}
\newcommand\TIinb{T^-_I}
\newcommand\Routb{R^+}
\newcommand\Routbsgnl{R^{\sign\lambda}}
\newcommand\Tout{T^+}
\newcommand\Toutb{T^+_{\ff}}
\newcommand\TIoutb{T^+_I}
\newcommand\Rlkf{(|\xib|^2-(\lambda-i0)^2)^{-1}}
\newcommand\Rlk{\rho_0(\lambda)}
\newcommand\Rmlk{\rho_0(-\lambda)}
\newcommand\Rpmlk{\rho_0(\pm\lambda)}
\newcommand\Rlka{\rho_1(\lambda)}
\newcommand\Rlkb{\rho_2(\lambda)}
\newcommand\Rilk{\rho_i(\lambda)}
\newcommand\reduced{\natural}
\newcommand\Rlf{R_0(\lambda)}
\newcommand\Rla{R_1(\lambda)}
\newcommand\Rlb{R_2(\lambda)}
\newcommand\Ril{R_i(\lambda)}
\newcommand\Rlj{R_j(\lambda)}
\newcommand\Rlft{R_0(\lambda)}
\newcommand\Rflambda{R_0^{\reduced}(\sigma)}
\newcommand\RV{R^{\reduced}_V}
\newcommand\Rfsigma{R_0^{\reduced}(\sigma)}
\newcommand\Rfsigmah{R_0^{\reduced}(\sigma^{1/2})}
\newcommand\Rfzero{R_0^{\reduced}(0)}
\newcommand\RlV{R^{\reduced}_V(\sigma)}
\newcommand\RlVi{R^{\reduced}_{V_i}(\sigma)}
\newcommand\RlVt{R_V(\lambda)}
\newcommand\RlVtL{{R}_V^L(\lambda)}
\newcommand\RlVtR{{R}_V^R(\lambda)}
\newcommand\RlVit{{R}_{V_i}(\lambda)}
\newcommand\RlVta{{R}_V^{(1)}(\lambda)}
\newcommand\RlVtk{{R}_V^{(k)}(\lambda)}
\newcommand\RlVatV{{R}_{V_{\alpha}}(\lambda)V_{\alpha}}
\newcommand\RlVatVa{{R}_{V_{\alpha_1}}(\lambda)V_{\alpha_1}}
\newcommand\RlVatVb{{R}_{V_{\alpha_2}}(\lambda)V_{\alpha_2}}
\newcommand\RlVatVk{{R}_{V_{\alpha_k}}(\lambda)V_{\alpha_k}}
\newcommand\RlVatVkk{{R}_{V_{\alpha_{k+1}}}(\lambda)V_{\alpha_{k+1}}}
\newcommand\RlVaptV{{R}_{V_{\alpha'}}(\lambda)V_{\alpha'}}
\newcommand\RlVapptV{{R}_{V_{\alpha''}}(\lambda)V_{\alpha''}}
\newcommand\RlVajtV{{R}_{V_{\alpha_j}}(\lambda)V_{\alpha_j}}
\newcommand\RlVaktV{{R}_{V_{\alpha_k}}(\lambda)V_{\alpha_k}}
\newcommand\RlVakktV{{R}_{V_{\alpha_{k+1}}}(\lambda)V_{\alpha_{k+1}}}
\newcommand\Tl{T(\lambda)}
\newcommand\Tlt{\tilde\Tl}
\newcommand\Tltp{\tilde T'(\lambda)}
\newcommand\Tltpp{\tilde T''(\lambda)}
\newcommand\Tli{T_i(\lambda)}
\newcommand\Tlit{\tilde\Tli}
\newcommand\Tlip{T_i'(\lambda)}
\newcommand\Tlipp{T_i''(\lambda)}
\newcommand\Tlj{T_j(\lambda)}
\newcommand\Tla{T_{\alpha}(\lambda)}
\newcommand\Tlaa{T_{\alpha_1}(\lambda)}
\newcommand\Tlab{T_{\alpha_2}(\lambda)}
\newcommand\Tlak{T_{\alpha_k}(\lambda)}
\newcommand\Tlakt{\tilde\Tlak}
\newcommand\Tlaj{T_{\alpha_j}(\lambda)}
\newcommand\Tlajj{T_{\alpha_{j+1}}(\lambda)}
\newcommand\Tlajp{T_{\alpha_j}'(\lambda)}
\newcommand\Tlajpt{\tilde\Tlajp}
\newcommand\Tlajt{\tilde\Tlaj}
\newcommand\Tlakk{T_{\alpha_{k+1}}(\lambda)}
\newcommand\Tlakkp{T_{\alpha_{k+1}}'(\lambda)}
\newcommand\Tlap{T_{\alpha'}(\lambda)}
\newcommand\Tlapt{\tilde\Tlap}
\newcommand\Tlapp{T_{\alpha''}(\lambda)}
\newcommand\Tkl{T^{(k)}(\lambda)}
\newcommand\Tcl{T^{\flat}(\lambda)}
\newcommand\Fl{F(\lambda)}
\newcommand\BlVt{\tilde B_V(\lambda)}
\newcommand\KBlVt{K_{\BlVt}}
\newcommand\BlVaat{B_{V_{\alpha_1}}(\lambda)}
\newcommand\BV{B_V}
\newcommand\Bone{B_1}
\newcommand\Btwo{B_2}
\newcommand\Bthree{B_3}
\newcommand\Banyj{B_j}
\newcommand\PlV{P_V(\lambda)}
\newcommand\PlVc{P_V^{\flat}(\lambda)}
\newcommand\Pl{P_0(\lambda)}
\newcommand\SVl{S_V(\lambda)}
\newcommand\Sjr{S_j^{\reduced}}
\newcommand\Rkp{{\mathcal R}^k_+}
\newcommand\Rkm{{\mathcal R}^k_-}
\newcommand\Rkpm{{\mathcal R}^k_{\pm}}
\newcommand\Phys{{\mathcal P}}
\newcommand\Pc{\overline{\mathcal P}}
\newcommand\pip{\pi^{\perp}}
\newcommand\pipa{\pi_\partial}
\newcommand\gammapa{\gamma_\partial}
\newcommand\pipah{\hat\pi_\partial}
\newcommand\pit{\tilde\pi}
\newcommand\xit{\tilde\xi}
\newcommand\zetat{\tilde\zeta}
\newcommand\etat{\tilde\eta}
\newcommand\sigmat{\tilde\sigma}
\newcommand\sigmahat{\hat\sigma}
\newcommand\thetat{\tilde\theta}
\newcommand\psit{\tilde\psi}
\newcommand\phit{\tilde\phi}
\newcommand\chit{\tilde\chi}
\newcommand\rhot{\tilde\rho}
\newcommand\xib{\bar\xi}
\newcommand\zetab{\bar\zeta}
\newcommand\thetab{\bar\theta}
\newcommand\etab{\bar\eta}
\newcommand\iotal{\iota_{\lambda}}
\newcommand\rhoat{\rhot_{\alpha_1}}
\newcommand\Lambdat{\tilde\Lambda}
\newcommand\Lambdati{\tilde\Lambda^{\text{in}}}
\newcommand\Lambdato{\tilde\Lambda^{\text{out}}}
\newcommand\Lambdatp{\tilde\Lambda^{\text{prop}}}
\newcommand\Lambdai{\Lambda^{\text{in}}}
\newcommand\Lambdao{\Lambda^{\text{out}}}
\newcommand\poles{\Lambda'}
\newcommand\rpoles{\Lambda_p}
\newcommand\thresholds{\Lambda}
\newcommand\Vt{\tilde V}
\newcommand\It{\tilde I}
\newcommand\half{{\frac{1}{2}}}
\newcommand\sigmah{\sigma^{1/2}}
\newcommand\bX{\partial X}
\newcommand\bXb{\partial \Xb}
\newcommand\Deltabt{\tilde\Delta_0}
\newcommand\strip{\Omega_T}
\newcommand\Kf{K^{\flat}}
\newcommand\Gs{G^{\sharp}}
\newcommand\Gt{\tilde G}
\newcommand\Osc{\sci\Omega}
\newcommand\OSc{{}^\Scl\Omega}
\newcommand\Osch{\sci\Omega^{\half}}
\newcommand\Oscmh{\sci\Omega^{-\half}}
\newcommand\Isc{I_{sc}}
\newcommand\os{{\text{os}}}
\newcommand\Qzl{Q^0_{-\lambda}}
\newcommand\Lie{{\mathcal L}}
\newcommand\bl{{\text b}}
\newcommand\scl{{\text{sc}}}
\newcommand\sccl{{\text{scc}}}
\newcommand\Scl{{\text{Sc}}}
\newcommand\ScLl{{\text{Sc,L}}}
\newcommand\ScRl{{\text{Sc,R}}}
\newcommand\Sccl{{\text{Scc}}}
\newcommand\sus{{\text{sus}}}
\newcommand\ssl{{\text{ss}}}
\newcommand\XXb{X^2_\bl}
\newcommand\XXbt{\Xt^2_\bl}
\newcommand\XXsc{X^2_\scl}
\newcommand\XXsct{\Xt^2_\scl}
\newcommand\XXSc{X^2_\Scl}
\newcommand\XXSct{\Xt^2_\Scl}
\newcommand\XXScL{X^2_\ScLl}
\newcommand\XXScR{X^2_\ScRl}
\newcommand\MMsc{M^2_\scl}
\newcommand\Deltab{\Delta_\bl}
\newcommand\Deltasc{\Delta_\scl}
\newcommand\DeltaSc{\Delta_\Scl}
\newcommand\DeltaScL{\Delta_\ScLl}
\newcommand\DeltaScR{\Delta_\ScRl}
\newcommand\prs{\sigma}
\newcommand\Nsc{N_\scl}
\newcommand\Nscp{N_{\scl,p}}
\newcommand\Nff{N_{\ff}}
\newcommand\Nffz{N_{\ff,0}}
\newcommand\Nffzp{N_{\ff,0,p}}
\newcommand\Nffl{N_{\ff,l}}
\newcommand\Nffml{N_{\ff,-l}}
\newcommand\Nmf{N_{\mf}}
\newcommand\Nmfz{N_{\mf,0}}
\newcommand\Nmfl{N_{\mf,l}}
\newcommand\Nmfml{N_{\mf,-l}}
\newcommand\ffb{\operatorname{bf}}
\newcommand\Ffb{\operatorname{bf'}}
\newcommand\ffsc{\operatorname{sf}}
\newcommand\ffSc{\operatorname{sf_C}}
\newcommand\Ffsc{\operatorname{sf'}}
\newcommand\rff{\rho_{\ff}}
\newcommand\rmf{\rho_{\mf}}
\newcommand\rffb{\rho_{\ffb}}
\newcommand\rffsc{\rho_{\ffsc}}
\newcommand\rFfsc{\rho_{\Ffsc}}
\newcommand\rffSc{\rho_{\ffSc}}
\newcommand\rinf{\rho_{\infty}}
\newcommand\CL{C_L}
\newcommand\CR{C_R}
\newcommand\betab{\beta_\bl}
\newcommand\betasc{\beta_\scl}
\newcommand\betaSc{\beta_\Scl}
\newcommand\BetaSc{\bar\beta_\Scl}
\newcommand\betaScL{\beta_\ScLl}
\newcommand\betaScR{\beta_\ScRl}
\newcommand\ScT{{}^\Scl T^*}
\newcommand\SccT{{}^\Scl \bar T^*}
\newcommand\ScS{{}^\Scl S^*}
\newcommand\Tb{{}^\bl T}
\newcommand\Tsc{{}^\scl T}
\newcommand\TSc{{}^\Scl T}
\newcommand\CSc{C_\Scl}
\newcommand\Lambdasc{{}^\scl\Lambda}
\newcommand\XXXb{X^3_\bl}
\newcommand\XXXsc{X^3_\scl}
\newcommand\XXXSc{X^3_\Scl}
\newcommand\XXXScO{X^3_{\Scl,O}}
\newcommand\XXXScF{X^3_{\Scl,F}}
\newcommand\XXXScS{X^3_{\Scl,S}}
\newcommand\XXXScC{X^3_{\Scl,C}}
\newcommand\KDsc{\operatorname{KD^{\half}_\scl}}
\newcommand\KDSc{\operatorname{KD^{\half}_\Scl}}
\newcommand\KDScEF{\operatorname{KD^{E,F}_\Scl}}
\newcommand\Oh{\operatorname{\Omega^{\half}}}
\newcommand\WFSc{\WF_\Scl}
\newcommand\WFtSc{\WF_{\text 3sc}}
\newcommand\WFScmf{\WF_{\Scl,\mf}}
\newcommand\WFScff{\WF_{\Scl,\ff}}
\newcommand\WFScs{\WF_{\Scl,\prs}}
\newcommand\WFScp{\WF'_\Scl}
\newcommand\WFScmfp{\WF'_{\Scl,\mf}}
\newcommand\WFScffp{\WF'_{\Scl,\ff}}
\newcommand\WFScsp{\WF'_{\Scl,\prs}}
\newcommand\Diffscc{\Diff_\sccl}
\newcommand\DiffSc{\Diff_\Scl}
\newcommand\DiffScc{\Diff_\Sccl}
\newcommand\DiffscI{\Diff_{\scl,\text{I}}}
\newcommand\VscI{\Vf_{\scl,\text{I}}}
\newcommand\DiffsV{\operatorname{Diff}_{\sus(V)}}
\newcommand\DiffsVsc{\operatorname{Diff}_{\sus(V),\scl}}
\newcommand\DiffsVCsc{\operatorname{Diff}_{\sus(V)-C,\scl}}   
\newcommand\Psisc{\Psop_\scl}
\newcommand\Psiscc{\Psop_\sccl}
\newcommand\Psiss{\Psop_\ssl}
\newcommand\Psisch{\Psop_{\scl,h}}
\newcommand\Psiscch{\Psop_{\sccl,h}}
\newcommand\PsiSc{\Psop_\Scl}
\newcommand\PsiScph{\Psop_{\Scl,\phi}}
\newcommand\PsiScra{\Psop_{\Scl,\rho^\sharp_a}}
\newcommand\PsiScc{\Psop_\Sccl}
\newcommand\PsiSccml{\Psop^{m,l}_\Sccl}
\newcommand\PsiScxx{\Psop^{*,*}_\Scl}
\newcommand\PsiScml{\Psop^{m,l}_\Scl}
\newcommand\PsiScmz{\Psop^{m,0}_\Scl}
\newcommand\PsiScmmz{\Psop^{-m,0}_\Scl}
\newcommand\PsiSckz{\Psop^{k,0}_\Scl}
\newcommand\PsiScmmml{\Psop^{-m,-l}_\Scl}
\newcommand\Psiscmkk{\Psop^{-k,k}_\scl}
\newcommand\Psiscmmmkk{\Psop^{-m-k,k}_\scl}
\newcommand\Psiscmoo{\Psop^{-1,1}_\scl}
\newcommand\Psiscmz{\Psop^{m,0}_\scl}
\newcommand\Psiscmmz{\Psop^{-m,0}_\scl}
\newcommand\PsiSckmkl{\Psop^{km,kl}_\Scl}
\newcommand\PsiScmplp{\Psop^{m',l'}_\Scl}
\newcommand\PsiScmmpllp{\Psop^{m+m',l+l'}_\Scl}
\newcommand\Psiscml{\Psop^{m,l}_\scl}
\newcommand\PsiScid{\Psop^{0,0}_\Scl}
\newcommand\PsiSczo{\Psop^{0,1}_\Scl}
\newcommand\PsiScmii{\Psop^{-\infty,\infty}_\Scl}
\newcommand\PsiScmiz{\Psop^{-\infty,0}_\Scl}
\newcommand\PsiScmoo{\Psop^{-1,1}_\Scl}
\newcommand\PsisCid{\Psop^{0,0}_{\scl-C}}
\newcommand\PsisC{\Psop_{\scl-C}}
\newcommand\Psiinf{\Psop_{\infty}}
\newcommand\Psiinfid{\Psop_{\infty}^0}
\newcommand\PsiFinf{\Psop_{\infty-\Fr}}
\newcommand\PsisVscml{\Psop^{m,l}_{\sus(V),\scl}}
\newcommand\PsisVsc{\Psop_{\sus(V),\scl}}
\newcommand\PsisVpsc{\Psop_{\sus(V_p),\scl}}
\newcommand\PsisVCSc{\Psop_{\sus(V)-C,\scl}}
\newcommand\SFinf{S_{\infty-\Fr}}
\newcommand\YsVC{Y^2_{\sus(V)-C,\scl}}
\newcommand\ffYsc{\ffsc_{\sus(V)}}
\newcommand\SXC{S(X;C)}
\newcommand\Ios{I_{\text{os}}}
\newcommand\pbL{\pi^2_{\bl,{\text L}}}
\newcommand\pbR{\pi^2_{\bl,{\text R}}}
\newcommand\pscL{\pi^2_{\scl,{\text L}}}
\newcommand\pscR{\pi^2_{\scl,{\text R}}}
\newcommand\PbO{\pi^3_{\bl,{\text O}}}
\newcommand\PscO{\pi^3_{\scl,{\text O}}}
\newcommand\PScO{\pi^3_{\Scl,{\text O}}}
\newcommand\PScF{\pi^3_{\Scl,{\text F}}}
\newcommand\PScC{\pi^3_{\Scl,{\text C}}}
\newcommand\PScS{\pi^3_{\Scl,{\text S}}}
\newcommand\pScL{\pi^2_{\Scl,{\text L}}}
\newcommand\pScR{\pi^2_{\Scl,{\text R}}}
\newcommand\CLF{\CL^F}
\newcommand\CLO{\CL^O}
\newcommand\CLS{\CL^S}
\newcommand\CLC{\CL^C}
\newcommand\DeltaYb{\Delta_{\bl,Y}}
\newcommand\DeltaYsc{\Delta_{\sus-\scl}}
\newcommand\diag{\operatorname{diag}}
\newcommand\Vf{{\mathcal V}}
\newcommand\Vb{{\mathcal V}_{\bl}}
\newcommand\Vsc{{\mathcal V}_{\scl}}
\newcommand\VSc{{\mathcal V}_{\Scl}}
\newcommand\VfI{\Vf_{\text{I}}}
\newcommand\VfIq{\Vf_{\text{I},q}}
\newcommand\scH{{}^\scl H}
\newcommand\scHg{\scH_g}
\newcommand\Hss{H_\ssl}
\newcommand\xh{\hat x}
\newcommand\Yh{\hat Y}
\newcommand\Zh{\hat Z}
\newcommand\Yb{\bar Y}
\newcommand\hb{\bar h}
\newcommand\xih{\hat\xi}
\newcommand\etah{\hat\eta}
\newcommand\muh{\hat\mu}
\newcommand\mub{\bar\mu}
\newcommand\nub{\bar\nu}
\newcommand\mubh{\widehat{\bar\mu}}
\newcommand\yb{\bar y}
\newcommand\zb{\bar z}
\newcommand\ub{\bar u}
\newcommand\Qb{\bar Q}
\newcommand\Wbp{{\bar W}^\perp}
\newcommand\Wp{W^\perp}
\newcommand\Kt{\tilde K}
\newcommand\Wt{\tilde W}
\newcommand\Ut{\tilde U}
\newcommand\yt{\tilde y}
\newcommand\ut{\tilde u}
\newcommand\vt{\tilde v}
\newcommand\ft{\tilde f}
\newcommand\htil{\tilde h}
\newcommand\St{\tilde S}
\newcommand\Pt{\tilde P}
\newcommand\Rt{\tilde R}
\newcommand\qt{\tilde q}
\newcommand\Qt{\tilde Q}
\newcommand\Xb{\bar X}
\newcommand\lambdat{\tilde\lambda}
\newcommand\betat{\tilde\beta}
\newcommand\epst{\tilde\epsilon}
\newcommand\ep{\epsilon}
\newcommand\bt{\tilde b}
\newcommand\Xt{\tilde X}
\newcommand\Mt{\tilde M}
\newcommand\At{\tilde A}
\newcommand\Et{\tilde E}
\newcommand\Ht{\tilde H}
\newcommand\at{\tilde a}
\newcommand\Ct{\tilde C}
\newcommand\pih{\hat\pi}
\newcommand\Rh{\hat R}
\newcommand\Ah{\hat A}
\newcommand\Bh{\hat B}
\newcommand\Ch{\hat C}
\newcommand\Gh{\hat G}
\newcommand\Hh{\hat H}
\newcommand\Qh{\hat Q}
\newcommand\Ph{\hat P}
\newcommand\Nh{\hat N}
\newcommand\Sh{\hat S}
\newcommand\Gcal{{\mathcal G}}
\newcommand\GcalC{{\mathcal G}_C}
\newcommand\Jcal{{\mathcal J}}
\newcommand\JcalC{{\mathcal J}_C}
\setcounter{secnumdepth}{3}
\newtheorem{lemma}{Lemma}[section]
\newtheorem{prop}[lemma]{Proposition}
\newtheorem{thm}[lemma]{Theorem}
\newtheorem{cor}[lemma]{Corollary}
\newtheorem{result}[lemma]{Result}
\newtheorem*{thm*}{Theorem}
\newtheorem*{prop*}{Proposition}
\newtheorem*{cor*}{Corollary}
\newtheorem*{conj*}{Conjecture}
\numberwithin{equation}{section}
\theoremstyle{remark}
\newtheorem{rem}[lemma]{Remark}
\newtheorem*{rem*}{Remark}
\theoremstyle{definition}
\newtheorem{Def}[lemma]{Definition}
\newtheorem*{Def*}{Definition}
\def\signature#1#2{\par\noindent#1\dotfill\null\\*
{\raggedleft #2\par}}

\renewcommand{\theenumi}{\roman{enumi}}
\renewcommand{\labelenumi}{(\theenumi)}

\title[Exponential decay]
{Exponential decay of eigenfunctions in many-body type scattering with
second order perturbations}
\author[Andras Vasy]{Andr\'as Vasy}
\address{Department of Mathematics, Massachusetts Institute of
Technology, MA 02139}
\email{andras@math.mit.edu}
\date{April 23, 2002}
\thanks{Partially supported by NSF grant \#DMS-99-70607.}
\
\begin{abstract}
We show the exponential decay of eigenfunctions of second-order
geometric many-body type
Hamiltonians at non-threshold energies. Moreover, in the case of first
order and small second order perturbations
we show that there are no eigenfunctions with positive
energy.
\end{abstract}
\maketitle

\section{Introduction and results}

In this paper we show that $L^2$-eigenfunctions of
elliptic second order many-body type perturbations $H$ of the Laplacian
with non-threshold eigenvalues $\lambda$ decay exponentially at a rate given by
the distance of $\lambda$ to the next threshold above it. If there are
no positive thresholds, this implies the super-exponential decay of
eigenfunctions at positive energies. We also show a unique continuation
theorem at infinity, namely that for first-order and small second-order
perturbations of the Laplacian, super-exponential decay of an eigenfunction
$\psi$ implies that $\psi$ is identically $0$. In particular, for these
perturbations, an inductive argument shows that such Hamiltonians have
no positive eigenvalues.
These generalize results of \cite{Vasy:Propagation-2}, where only potential
scattering was considered, although already in a geometric setting,
and the pioneering work of Froese and Herbst \cite{FroExp}, in which
they considered many-body potential scattering in Euclidean space.

The methods are closely related to both those of Froese and Herbst
and of the two-body type unique continuation theorems discussed in
\cite{Hormander:Uniqueness} and \cite[Theorem~17.2.8]{Hor}. However,
the geometric nature of the problem forces a systematic treatment
of various `error terms', and in particular the use of a very stable
argument. In particular, we emphasize throughout that
for the exponential decay results only the indicial operators of $H$,
which are non-commutative analogues of the usual principal symbol,
affect the arguments, hence $H$ can be generalized a great deal more.
In addition, for unique continuation result only the indicial operators of $H$
and its symbol in a high-energy sense (as in
`ellipticity with a parameter', or after rescaling, as in semiclassical
problems) play a role. This explains, in particular, the first order
(or small second-order)
hypothesis on the perturbations for the unique continuation theorem, and
raises the question to whether this theorem also holds under non-trapping
conditions on the metric near infinity, or even more generally. The
key estimates arise from a positive commutator estimate for the
conjugated Hamiltonian, which is closely related to H\"ormander's
solvability condition for PDE's \cite{Hormander:Differential,
Hormander:Solutions, Duistermaat-Sjostrand:Global}; see
\cite{Zworski:Numerical} for a recent
discussion, including the relationship to numerical computation.

Before stating the results precisely,
recall from \cite{RBMSpec}
that if $\Xb$ is a manifold with boundary and $x$ is
a boundary defining function on $\Xb$, a scattering metric $g_0$
is a Riemannian metric on $X=\Xb^\circ$ which is of the form
$g_0=x^{-4}\,dx^2+x^{-2}h$ near $\pa \Xb$, where $h$ is a symmetric
2-cotensor that restricts to a metric on $\pa \Xb$. Let $\Delta_{g_0}$
be the Laplacian of this metric. This is a typical element of
$\Diffsc(\Xb)$, the algebra of scattering differential operators.
The latter is generated, over $\Cinf(\Xb)$, by the vector fields
$\Vsc(\Xb)=x\Vb(\Xb)$; $\Vb(\Xb)$ being the Lie algebra of $\Cinf$ vector
fields on $\Xb$ that are tangent to $\pa \Xb$.

In this paper we consider many-body type Hamiltonians. That is, let
$\calC$ be a cleanly intersecting family of closed embedded
submanifolds of $\pa \Xb$
which is closed under intersections and which includes $C_0=\pa \Xb$ (the
latter only for convenient notation). As shown in \cite{Vasy:Propagation-Many},
one can resolve $\calC$
by blowing these up inductively, starting with the submanifold of the
lowest dimension. The resulting space $[\Xb;\calC]$ is a manifold with corners,
and the blow-down map $\beta:[\Xb;\calC]\to \Xb$ is smooth.
Then $\DiffSc(\Xb,\calC)$ is similar to $\Diffsc(\Xb)$, but with
coefficients that are in $\Cinf([\Xb;\calC])$: $\DiffSc(\Xb,\calC)=
\Cinf([\Xb;\calC])\otimes_{\Cinf(\Xb)}\Diffsc(\Xb)$. More generally, if $E,F$
are vector bundles over $\Xb$, we can consider differential operators mapping
smooth sections of $\beta^*E$ to smooth sections of $\beta^*F$, denoted by
$\DiffSc(\Xb,\calC;E,F)$, or simply $\DiffSc(\Xb,\calC;E)$ if $E=F$.
The vector fields in $\DiffSc(\Xb;\calC)$ form exactly the set of all
smooth sections of a vector bundle, denoted by $\TSc[\Xb;\calC]$ over
$[\Xb;\calC]$, namely the pull-back of $\Tsc\Xb$ by $\beta$. The dual
bundle is denoted $\ScT[\Xb;\calC]$.

It may help the reader if we explain why the Euclidean setting is a particular
example of this setup. Indeed, the reader may be interested in the Euclidean
magnetic and metric scattering specifically; if so, all the arguments
given below can be translated into Euclidean notation as follows.
There $X$ is a vector space with a metric $g_0$,
which can hence by identified with $\Rn$. Moreover, $\Xb$ is the radial
(or geodesic) compactification of $\Rn$ to a ball. Explicitly, this
arises by considering `inverse' polar coordinates, and writing $w\in X$ as
$w=r\omega=x^{-1}\omega$, $\omega\in\Sn$, so $x=|w|^{-1}$,
e.g.\ in $|w|\geq 1$. In particular,
$\pa \Xb$ is given by $x=0$, i.e.\ it is just $\Sn$. The metric $g_0$
then has the form $dr^2+r^{-2}h_0=x^{-4}dx^2+x^2h_0$, where $h_0$ is the
standard metric on $\Sn$, so $(\Xb,g_0)$
fits exactly into this framework. Moreover,
in the many-body setting, one is given a collection $\calX=\{X_a:\ a\in I\}$
of linear subspaces of $X$. The corresponding cleanly intersecting family
is given by $\calC=\{C_a:\ a\in I\}$, where $C_a=\Xb_a\cap\pa\Xb$, and
$\Xb_a$ is the closure of $X_a$ in $\Xb$. Thus, $C_a$ can also be thought
of as the intersection of the unit sphere in $\Rn$ with $X_a$. Then
$\TSc[\Xb;\calC]$, $\ScT[\Xb;\calC]$ are trivial vector bundles over
$[\Xb;\calC]$; namely $\ScT[\Xb;\calC]=[\Xb;\calC]\times X^*$, $X^*$ being
the dual vector space of $X$.

We can now describe the operators $H$ we consider in this paper.
First, we assume that $H=\Delta\otimes\Id_E+V$ where
$\Delta=\Delta_g$ is the Laplacian
of a metric $g$ such that
\begin{equation*}
g\in\Cinf([\Xb;\calC];\ScT[\Xb;\calC]\otimes\ScT[\Xb;\calC])
\end{equation*}
is symmetric, $g-g_0$ vanishes at the free face, i.e.\ the lift of $C_0$,
and $V\in\DiffSc^1(\Xb;\calC;E)$ is formally self-adjoint and vanishes
at the free face.

Now, $g_0$ induces an orthogonal decomposition of
$\sct \Xb$ at each $C_a\in\calC$, which, with the Euclidean notation
corresponds to the decomposition
\begin{equation*}
T^*_{(w_a,w^a)}X_0=T^*_{w_a}X_a\oplus
T^*_{w^a}X^a;
\end{equation*}
$X^a$ being the orthocomplement of $X_a$.
In the geometric setting, $X^a$ is simply
short hand for
the fibers $\beta_a^{-1}(y_a)$, $y_a\in C_a$
of the front face of $\beta_a:[\Xb;C_a]\to\Xb$,
while a neighborhood of infinity in
the radial compactification $\Xb_a$ of
$X_a$ stands for $C_a\times[0,\ep)_x$, see \cite{Vasy:Propagation-Many}.
This allows us to define the indicial operators
of differential operators
$A\in\DiffSc^m(\Xb;\calC)$ invariantly (even in the geometric setting),
at a cluster $a$, as a family of operators $\hat A_a(y_a,\xi_a)$, $y_a\in C_a$,
$\xi_a\in \sct_{y_a}\Xb_a$, on functions on $X^a$, by
freezing the coefficients of $A$ at $y_a$ and replacing derivatives
$D_{(w_a)_j}$ by $(\xi_a)_j$. (Technically the $a$-indicial operators are
defined on a blow-up of $C_a$, i.e.\ the above definition is valid for
$y_a\in C_{a,\reg}$, i.e.\ for $y_a$ away from all $C_b$ which do not satisfy
$C_b\supset C_a$; see \cite{Vasy:Propagation-Many} for the detailed setup.)
For example, the $a$-indicial operator of $\Delta_{g_0}-\lambda$ is
$|\xi_a|^2_{y_a}+\Delta_{X^a}(y_a)-\lambda$; here $\Delta_{X^a}(y_a)$ is,
for each $y_a\in C_a$,
a translation invariant operator on the fibers of the front face of
$\beta_a:[\Xb;C_a]\to \Xb$.

We also assume that restricted to a front face, $g-g_0$ is a section
of $T^*X^a\otimes T^*X^a$, i.e.\ depends on the interaction variables
only over $C_a$, and that for each $y_a\in C_a$, $\hat V_a(y_a,\xi_a)
\in\Diff^1(X^a;E)$ is independent of $\xi_a$.
In other words, we assume that
all indicial operators of $H$ are pointwise in $C_a$
`product-type', i.e.\ have
the form
\begin{equation}\label{eq:H_a-ind-pt-pr}
\hat H_a(\xi_a)=|\xi_a|^2_{y_a}\otimes\Id_E+H^a(y_a),
\end{equation}
where $H^a(y_a)\in\Diff^2(X^a;E)$, and $y_a$ is the variable along $C_a$.

This product assumption is sufficient for all of our results, provided
that on the complement of the set of thresholds, a Mourre-type global
positive commutator estimate (local only in the spectrum of $H$) holds.
If $H^a(y_a)$ has $L^2$-eigenvalues, the existence of such a
Mourre estimate at various energies certainly depends on the behavior
of the eigenvalues as a function of $y_a$, as shown by a related problem
involving scattering by potentials of degree zero, \cite{Herbst:Spectral,
Herbst-Skibsted:Quantum, Hassell-Melrose-Vasy:Spectral}.
So we assume in this paper that either $H^a(y_a)$ has no $L^2$
eigenvalues for any $a$, or in a neighborhood of $C_a$, $\Xb$ has
the structure of (the radial compactification of)
a conic slice of $X_a\times X^a$, over which $E$ is
trivial, and the indicial operators satisfy
\begin{equation}\label{eq:ind-H_a-strong-pr}
\hat H_a(\xi_a)=|\xi_a|^2\otimes\Id_E+H^a,
\end{equation}
so $H^a$
is independent of $\xi_a$, and in particular of its projection $y_a$ to $C_a$.
Here $H^a$ is called the subsystem Hamiltonian for the subsystem $a$; it
is also a many-body Hamiltonian (but one corresponding to fewer particles
in actual many-body scattering!).

Examples of such Hamiltonians include the Laplacian of
metric perturbations of $g_0$ in the Euclidean setting, both on functions, and
more generally on forms, and the square
of associated self-adjoint Dirac operators.
Namely, let $X_a$ be the collision planes, and let
$g^a\in\Cinf_c(X^a;T^*X^a\otimes T^*X^a)$, $a\neq 0$, be symmetric.
Then $g^a$ can also be regarded as a section of $T^*X\times T^*X$,
and $g=g_0+\sum_a g^a$ satisfies these criteria provided that it is
positive definite (i.e.\ a metric). Of course, the
compact support of the $g^a$ can be replaced by first order
decay at infinity as a section of $\sct \Xb^a\otimes\sct \Xb^a$,
i.e.\ relative to the translation-invariant basis
$dw^a_j\otimes dw^a_k$, $j,k=1,\ldots\dim X^a$, of $\Cinf(X^a;T^*X^a
\otimes T^*X^a)$.

We prove the following results. Let $\Lambda$
denote the thresholds of $H$, i.e.\ the set of the $L^2$ eigenvalues
of all of its subsystem Hamiltonians $H^a$;
this is a closed countable subset of
$\Real$. (Note that $0$ is an eigenvalue of $H^0$, hence it is a threshold
of $H$; if no non-trivial subsystem has an $L^2$-eigenvalue then
$\Lambda=\{0\}$.)
The following theorem states that non-threshold eigenfunctions
of $H$ decay exponentially at a rate given by the distance of the eigenvalue
from the nearest threshold above it. This result also explains why
the unique continuation theorem is considered separately, namely why
super-exponential decay assumptions are natural there (unlike the Schwartz
assumptions of \cite[Theorem~17.2.8]{Hor} in the two-body type setting).

\begin{thm*}$[$cf.\ \cite[Proposition~B.2]{Vasy:Propagation-2} and
Froese and Herbst, \cite[Theorem 2.1]{FroExp}$]$
Let $\lambda\in\Real\setminus\Lambda$, and suppose that
$\psi\in L^2_{\scl}(\Xb)$ satisfies $H\psi=\lambda\psi$. Then $e^{\alpha/x}\psi
\in L^2_{\scl}(\Xb)$ for all $\alpha\in\Real$ such that
$[\lambda,\lambda+\alpha^2]\cap\Lambda=\emptyset$, i.e.
\begin{equation*}
\sup\{\lambda+\alpha^2:\ e^{\alpha/x}\psi\in L^2_{\scl}(\Xb)\}
\geq\inf\{\lambda'\in\Lambda:\ \lambda'>\lambda\}.
\end{equation*}
\end{thm*}

The estimates leading to this theorem are uniform, and in fact yield,
as observed by Perry \cite{Perry:Exponential}
in the Euclidean many-body potential scattering, that eigenvalues
cannot accummulate at thresholds from above, hence the following corollary.

\begin{cor*}
The thresholds $\lambda\in\Lambda$ are isolated from above, i.e.\ for
$\lambda\in\Lambda$ there exists $\lambda'>\lambda$ such that
$(\lambda,\lambda')\cap(\Lambda\cup\pspec(H))=\emptyset$.
\end{cor*}

The unique continuation theorem at infinity is the following.

\begin{thm*}$[$cf.\ \cite[Proposition~B.3]{Vasy:Propagation-2} and
Froese and Herbst, \cite[Theorem 3.1]{FroExp}$]$
Let $\lambda\in\Real$ and let $d$ denote a metric giving the usual topology on
$\Cinf$ sections of $\ScT[X;\calC]\otimes\ScT[X;\calC]$.
There exists $\ep>0$ such that if
$d(g,g_0)<\ep$ and $H\psi
=\lambda\psi$, $\exp(\alpha/x)\psi\in L^2_{\scl}(\Xb)$ for all $\alpha$,
then $\psi=0$.
\end{thm*}

As an immediate corollary we deduce the absence of positive eigenvalues
for first order perturbation and small second order perturbations
of $\Delta_{g_0}$.

\begin{thm*}
Let $\lambda>0$, $g$ is close to $g_0$ in a $\Cinf$ sense. Suppose that $H\psi
=\lambda\psi$, $\psi\in L^2_{\scl}(\Xb)$.
Then $\psi=0$.
\end{thm*}

\begin{proof}
One proceeds inductively, showing that $H^a$ does not have any positive
eigenvalues, starting with $a=0$, when this is certainly true. So suppose
that for all $b$ such that $X^b\subsetneq X^a$, $H^b$ does not have
any positive eigenvalues. Then the set of thresholds for $H^a$ is
disjoint from $(0,+\infty)$, so by the first theorem any eigenfunction
with a positive eigenvalue decays super-exponentially, and then by
the second theorem it vanishes. This completes the inductive step.
\end{proof}

\begin{rem*}
The last result in particular applies to $H=\Delta_g$ on functions
even if $g$ restricted to the
front face of $[\Xb;C_a]\to \Xb$ is any smooth section of $T^*X^a
\otimes T^*X^a$, i.e.\ only \eqref{eq:H_a-ind-pt-pr} holds (rather than
\eqref{eq:ind-H_a-strong-pr}) to show that $H$ has no $L^2$
eigenfunctions at all. Indeed, proceeding inductively
as in the proof, we may assume that
for all $b$ such that $C_a\subsetneq C_b$, $\hat H_b(\xi_b)$ does not have
any $L^2$-eigenvalues. Thus, the Mourre-type estimate is valid,
hence $\hat H_a(\xi_a)$ has no positive eigenvalues.
But $\hat H_a(\xi_a)\geq 0$, so it cannot have
negative energy bound states, and by elliptic regularity, any
$L^2$ element $\psi$ of its nullspace would be in $\Hsc^\infty(\Xb^a)$ thus
$\Delta_g=(d+\delta)^2=\delta d+d\delta$ shows that $d\psi=0$, hence
$\psi=0$. Thus, $\hat H_a(\xi_a)$ has no $L^2$-eigenvalues, completing
the inductive step.
\end{rem*}

The rough idea of the proof of the two main results
is to conjugate by exponential weights
$e^F$, where $F$ is a symbol of order $1$, for example $F=\alpha/x$.
If $\psi$ is an eigenfunction of $H$ of eigenvalue $\lambda$, then
$\psi_F=e^F\psi$ solves
\begin{equation*}
P\psi_F=0\Mwhere P=H(F)-\lambda=e^FHe^{-F}-\lambda.
\end{equation*}
Now $\re P$ is given by $H-\alpha^2-\lambda$, while $\im P$ is given by
$-2\alpha(x^2 D_x)$, modulo $x\DiffSc(\Xb,\calC)$. By elliptic regularity,
using $P\psi_F=0$, $\|\psi_F\|_{x^p\Hsc^k(\Xb)}$ is bounded by
$C_{k,p}\|\psi_F\|_{x^pL^2_{\scl}(\Xb)}$, so the order of various
differential operators can be neglected, while the weight is important.
Since
\begin{equation*}
P^*P=(\re P)^2+(\im P)^2+i(\re P\im P-\im P\re P),
\end{equation*}
so
\begin{equation}\label{eq:comm-8-0}
0=(\psi_F,P^*P\psi_F)=\|\re P\psi_F\|^2+\|\im P\psi_F\|^2+
(\psi_F,i[\re P,\im P]\psi_F).
\end{equation}
Now, $[\re P,\im P]\in x\Diffsc^2(\Xb,\calC)$, i.e.\ has an extra
order of vanishing, which shows that
\begin{equation*}
\|\re P\psi_F\|\leq C_1\|x^{1/2}\psi_F\|,
\ \|\im P\psi_F\|\leq C_1\|x^{1/2}\psi_F\|.
\end{equation*}
Due to the extra factor of $x^{1/2}$, this can be interpreted roughly as
$\psi_F$ being close to being in the nullspace of both $\re P$ and of
$\im P$, hence both of $H-\lambda-\alpha^2$ and $x^2D_x$.

If, moreover, $(\psi_F,i[\re P,\im P]\psi_F)$ is positive, modulo
terms involving $\re P$ and $\im P$ (which can be absorbed in the
squares in \eqref{eq:comm-8-0}), and terms of the form $(\psi_F,R\psi_F)$,
$R\in x^2\DiffSc(\Xb,\calC)$, which are thus bounded by $C_2\|x\psi_F\|^2$,
then the factor $x$ (which has an extra $x^{1/2}$ compared to
$\|x^{1/2}\psi_F\|$) yields easily a bound for $\|x^{1/2}\psi_F\|$ in
terms of $\|\psi\|$. This gives estimates for the
norm $\|x^{1/2}\psi_F\|$, uniform both in $F$ and in $\psi$. A regularization
argument in $F$ then gives the exponential decay of $\psi$.

The positivity of $(\psi_F,i[\re P,\im P]\psi_F)$, in the sense described
above, is easy to see if we replace $i[\re P,\im P]$ by
$i[H-\lambda-\alpha^2,-2\alpha x^2D_x]$: this commutator is a standard one
considered in many-body scattering, although the even more usual one
would be $i[H-\lambda-\alpha^2,-2xD_x]$, whose local positivity
in the spectrum of $H$ is the Mourre estimate
\cite{Mourre:Operateurs, Perry-Sigal-Simon:Spectral,
FroMourre}. Indeed, the latter
commutator is the one considered by Froese and Herbst in Euclidean
many-body potential scattering, and we could adapt their argument
(though we would need to deal with numerous error terms) to our setting.
However, the argument presented here is more robust, especially in the
high energy sense discussed below, in which their approach would not
work in the generality considered here. There is one exception: for
$\alpha=0$, $\im P$ degenerates, and in this case we need to `rescale'
the commutator argument, and consider $i[H-\lambda-\alpha^2,-2xD_x]$
directly.

We next want to let
$\alpha\to\infty$. Since most of the related literature considers
semiclassical problems, we let $h=\alpha^{-1}$, and replace $P$ above
by $P_h=h^2P$, which is a semiclassical differential operator, $P_h
\in\Diff_{\Scl,h}^2(\Xb,\calC)$. Here
\begin{equation*}
\Diff_{\Scl,h}(\Xb,\calC)
=\Cinf([\Xb;\calC])\otimes_{\Cinf(\Xb)}\Diff_{\scl,h}(\Xb),
\end{equation*}
and $\Diff_{\scl,h}(\Xb)$ is the algebra of semiclassical scattering
differential operators discussed, for example, in
\cite{Vasy-Zworski:Semiclassical}. It is generated by $h\Vsc(\Xb)$ over
$\Cinf(\Xb\times[0,1)_h)$.
In this semiclassical sense, the first and
zeroth order terms in $H$ do not play a role in $P_h$: their
contribution is in $h\Diff_{\Scl,h}^1(\Xb,\calC)$, hence their contribution
to the commutator $i[\re P_h,\im P_h]$ is in $xh^2\Diff_{\Scl,h}(\Xb,\calC)$.
Moreover, if $g$ is close to $g_0$, then $i[\re P_h,\im P_h]$ is close
to the corresponding commutator with $P_h$ replaced by $h^2(e^F \Delta_{g_0}
e^{-F}-\lambda)$. Since in the latter case the commutator is
positive, modulo terms than can be absorbed in the
two squares in \eqref{eq:comm-8-0}, $i[\re P_h,\im P_h]$ is also positive
for $g$ near $g_0$. This gives an estimate as above, from which
the vanishing of $\psi$ near $x=0$ follows easily.

We remark that the estimates we use are related
to the usual proof of unique continuation at infinity on $\Rn$ (i.e.\ not
in the many-body setting), see \cite[Theorem~17.2.8]{Hor}, and to
H\"ormander's solvability condition for PDE's in terms of the
real and imaginary parts of the principal symbol.
Indeed, although
in \cite[Theorem~17.2.8]{Hor} various changes of coordinates are used first,
which change the nature of the PDE at infinity,
ultimately the necessary estimates also arise from a commutator of the
kind $i[\re P,\im P]$. However, even in that setting,
the proof we present appears more natural
from the point of view of scattering than the one presented there,
which is motivated by unique continuation at points in $\Real^n$.
In particular, the reader who is interested in the setting of
\cite[Theorem~17.2.8]{Hor} should be able to skip the proof of
Theorem~\ref{thm:Mourre-est}, which is rather simple (a Poisson bracket
computation) in that case. We remark
that related estimates, obtained by different techniques, form the
backbone of the (two-body type)
unique continuation results of Jerison and Kenig
\cite{Jerison-Kenig:Unique, Jerison:Carleman}.

The true flavor of our arguments is most
clear in the proof of the unique continuation theorem,
Theorem~\ref{thm:unique}. The reason is that on the one hand there is
no need for regularization of $F$, since we are assuming super-exponential
decay, on the other hand the positivity of $i[\re P_h,\im P_h]$ is
easy to see.

The structure of the paper is the following. In Section~2 we discuss
various preliminaries, including the structure of the conjugated
Hamiltonian and a Mourre-type global positive commutator estimate.
In Section~3 we prove the exponential decay of non-threshold eigenfunctions.
In Section~4, we prove the unique continuation theorem at
infinity. Finally, for the sake of completeness, and since technically
the usual statements of the Mourre estimate do not discuss the present setting,
we include its proof in the Appendix.
We emphasize that the presence of bundles such as $E$ makes
no difference in the discussion, hence they are ignored in order to
keep the notation manageable; see Remarks~\ref{rem:FH} and \ref{rem:Mourre}
for further information.

I am very grateful to Rafe Mazzeo, Richard Melrose, Daniel Tataru and Maciej
Zworski for helpful discussions. I also thank Rafe Mazzeo for a careful reading
of the manuscript, and his comments which improved it significantly.

\section{Preliminaries}

We first remark that
the Riemannian density of a metric $g$ has the form
\begin{equation}\label{eq:dg-form}
dg=\sqrt{\det (g_{ij})}\,dx\,dy=\tilde g\,\frac{dx\,dy}{x^{n+1}}\,\ n=\dim X,
\ \tilde g\in\Cinf([X;\calC]).
\end{equation}
By our conditions on the form of $g$, the Laplacian takes the following form
\begin{equation*}
\Delta_g=(x^2 D_x)^2+\sum_j b_j x^2 P_j+\sum_j x c_j Q_j+xR
\end{equation*}
with $b_j, c_j\in\Cinf([\Xb;\calC])$,
$P_j\in\Diff^2(\bXb)$, $Q_j\in\Diff^1(\bXb)$, $R\in \DiffSc^2(\Xb,\calC)$.
Hence, $H=\Delta_g+V$ takes the form
\begin{equation*}
H=(x^2 D_x)^2+\sum_j b'_j x^2 P'_j+\sum_j x c'_j Q'_j+e+xR',
\end{equation*}
with $b'_j, c'_j,e\in\Cinf([\Xb;\calC])$,
$P'_j\in\Diff^2(\bXb)$, $Q'_j\in\Diff^1(\bXb)$, $R'\in \DiffSc^2(\Xb,\calC)$.

Below we consider the conjugated Hamiltonian $H(F)=e^{F}He^{-F}$, where
$F$ is a symbol of order $1$. The exponential weights will facilitate
exponential decay estimates, and eventually the proof of unique
continuation at infinity.
Let $x_0=\sup_{\Xb} x$. By altering $x$ in a compact subset of $X$, we
may assume that $x_0<1/2$; we do this for the convenience of notation below.
We let
$S^m([0,1)_x)$ is the space of all symbols $F$ of order $m$ on $[0,1)$,
which satisfy $F\in\Cinf((0,1))$, vanish on $(1/2,1)$, and for which
$\sup|x^{m+k}\partial_x^k F|<\infty$ for all $k$. The topology of $S^m
([0,1))$ is given by the seminorms $\sup|x^{m+k}\partial_x^k F|$. Also,
the spaces $S^m(\Xb)$, resp.\ $S^m([\Xb;\calC])$, of symbols is defined
similarly, i.e.\ it is
given by seminorms $\sup|x^m P F|$, $P\in\Diffb^k(\Xb)$, resp.\ $P
\in\Diffb^k([\Xb;\calC])$. In the following
lemma $\Diffscc(\Xb)$, as usual, stands for non-classical (non-polyhomogeneous)
scattering differential operators (i.e.\ scattering differential operators
with non-polyhomogeneous coefficients), corresponding to the lack of
polyhomogeneity of $F$. In particular, $\Diffscc^0(\Xb)=S^0(\Xb)$ (considered
as multiplication operators). Similarly, $\DiffScc(\Xb,\calC)
=S^0([\Xb;\calC])\otimes_{S^0(\Xb)}\Diffscc(\Xb)$ stands for the
corresponding calculus of many-body differential operators.

\begin{lemma}\label{lemma:FH-1}
Suppose $\lambda\in\Real$,
$H\psi=\lambda\psi$, $\psi\in L^2_{\scl}(\Xb)$. Suppose also
that $\alpha\geq 0$, and for all $\beta$ we have $x^{-\beta}\exp(\alpha/x)
\psi\in L^2_{\scl}(\Xb)$. Then with $F\in S^1([0,1))$,
$F\leq\alpha/x+\beta|\log x|$
for some $\beta$, $\supp F\subset [0,1/2)$,
$\psi_F=e^F\psi=e^{F(x)}\psi$,
\begin{equation*}
P=P(F)=e^F(H-\lambda)e^{-F}=H(F)-\lambda,\ H(F)=H+e^F[H,e^{-F}],
\end{equation*}
we have $\psi_F\in\dCinf(\Xb)$,
\begin{equation}\label{eq:FH-1}
P(F)\psi_F=0,
\end{equation}
\begin{equation}\label{eq:FH-1a}
P(F)=H-2(x^2D_x F)(x^2D_x)+(x^2D_x F)^2-\lambda+xR_1,\quad R_1\in
\DiffScc^2(\Xb,\calC),
\end{equation}
with
\begin{equation}\label{eq:FH-1b}
\re P(F)=H+(x^2D_x F)^2-\lambda+xR_2,\ \im P(F)=2(x^2\pa_x F)(x^2D_x)+xR_3,
\end{equation}
$R_2,R_3\in\DiffScc^2(\Xb,\calC)$, $R_j$ bounded as long as
$x^2\pa_x F$ is bounded in $S^0([0,1))$, hence as long as
$F$ is bounded in $S^1([0,1))$.
The coefficients of the $xR_2$, $xR_3$ are in fact polynomials with vanishing
constant term, in $(x^2\pa_x)^{m+1} F$, $m\geq 0$.

In particular,
\begin{equation}\begin{split}\label{eq:FH-3}
i[\re P(F),\im P(F)]=&i[H+(x^2D_x F)^2,2(x^2\pa_x F)(x^2D_x)]\\
&+\re P(F) xR_4+xR_5\re P(F)\\
&+\im P(F) xR_6+xR_7\im P(F)+x^2R_8,
\end{split}\end{equation}
where $R_j\in \DiffScc^2(\Xb,\calC)$ are bounded as long as $x^2\pa_x F$
is bounded
in $S^0([0,1))$.
\end{lemma}

\begin{rem}\label{rem:FH}
All but the first term on the right hand side of \eqref{eq:FH-3}
should be considered error terms, even though they are only of the
same order (in terms of decay at $\pa\Xb$) as
\begin{equation*}
i[H+(x^2D_x F)^2-\lambda,2(x^2\pa_x F)(x^2D_x)],
\end{equation*}
due to the lack of commutativity of $\DiffScc(\Xb,\calC)$ even to top order.
The reason is that these terms contain factors of $\re P(F)$ and
$\im P(F)$, and we will have good control over $\re P(F)\psi_F$
and $\im P(F)\psi_F$.

For similar reasons, the presence of bundles $E$
 would make no difference, since
even if they are present, \eqref{eq:FH-1b} is unaffected, hence
\eqref{eq:FH-3} holds as well.
\end{rem}

\begin{proof}
First note that
\begin{equation}\label{eq:FH-p1}
[x^2D_x, e^F]=(x^2 D_x F)e^F,\qquad x^2 D_x F\in S^0([0,1)),
\end{equation}
so $e^F[H,e^{-F}]\in\Diffscc^1(\Xb,\calC)$. The dependence of the terms of
$P(F)$ on $F$ thus comes from $x^2 D_x F$, and its commutators through
commuting it through other vector fields (as in rewriting
$(x^2D_x)(x^2 D_x F)$ as $(x^2 D_x F)(x^2D_x)$ plus a commutator term), hence
through $(x^2 D_x)^{m+1} F$, $m\geq 0$.

Now, \eqref{eq:FH-1}, which a priori holds in a
distributional sense, $\psi_F\in x^rL^2_\scl(\Xb)$ for all $r$,
and the ellipticity of $\prs_{\Scl,2}(H)$ show that
$\psi_F\in\dCinf(\Xb)$.

We use
\begin{equation*}\begin{split}
\re P(F)=\half(P(F)+P(F)^*)&=H-\lambda+\half (e^F[H,e^{-F}]-[H,e^{-F}]e^F)\\
&=H-\lambda+\half[e^F,[H,e^{-F}]]
\end{split}\end{equation*}
to prove \eqref{eq:FH-1b} (note that only the $(x^2D_x)^2$ terms in $H$ gives
a non-vanishing contribution to the double commutator).
Finally, \eqref{eq:FH-3} follows since
\begin{equation*}
Q=\re P(F)-(H+(x^2D_x F)^2-\lambda)\in 
x\DiffSc^2(\Xb,\calC),
\end{equation*}
so $[Q,\im P(F)]=Q\im P(F)-\im P(F) Q$ is of the
form of the $R_4$ and $R_5$ terms, and similarly for
$Q'=\im P(F)-2(x^2\pa_x F)(x^2D_x)$.
\end{proof}

In light of \eqref{eq:FH-3}, we need a positivity result for $i[x^2D_x,H]$.
Such a result follows directly from a Poisson bracket computation if
$H$ is a geometric 2-body type operator. In general, it requires a
positive commutator estimate that is closely related to, and can
be readily deduced from, the well-known
Mourre estimate \cite{Mourre:Operateurs, Perry-Sigal-Simon:Spectral,
FroMourre},
whose proof goes through in this generality. We will
briefly sketch its proof in the appendix for the sake of completeness.
So let
\begin{equation*}
d(\lambda)=\inf\{\lambda-\lambda':\lambda'\leq\lambda,\ \lambda'\in\Lambda\}
\end{equation*}
be the distance of $\lambda$ to the next threshold below it.
Let $\chi\in\Cinf_c([0,1))$ be supported near $0$, identically $1$ on
a smaller neighborhood of $0$, and let
\begin{equation*}
B=\half(\chi(x) x^2D_x+(\chi(x) x^2D_x)^*)
\end{equation*}
be the symmetrization of the radial vector field.
Now $x^2D_x$ is formally self-adjoint with respect to the measure
$\frac{dx\,dy}{x^2}$,
and if $C$ is formally
self-adjoint with respect to a density $dg'$ then its adjoint
with respect to $\alpha\,dg'$, $\alpha$ smooth real-valued,
is $\alpha^{-1}C\alpha=C+\alpha^{-1}[C,\alpha]$.
Since $xD_x$ is tangent to all elements of $\calC$,
$xD_x (x^{-n+1}\tilde g)\in x^{-n+1}\Cinf([\Xb;\calC])$.
In particular,
\begin{equation*}
ib=x^{-1}(B-\chi(x) x^2D_x)\in \Cinf([\Xb;\calC])
\end{equation*}
and $b$ is real-valued. It is easy to check that
in the particular case when $g=g_0$ then
$b\in\Cinf(\Xb)$ and $b|_{\pa X}=\frac{n-1}{2}$, $n=\dim X$.

The first order differential operator $B$ gives rise to:

\begin{thm}\label{thm:Mourre-est}
Suppose $\lambda\in\Real\setminus\Lambda$. For $\ep>0$ there exists $\delta>0$
such that for $\phi\in\Cinf_c(\Real;[0,1])$ with
$\supp\phi\subset(\lambda-\delta,
\lambda+\delta)$,
\begin{equation*}
\phi(H)i[B,H]\phi(H)\geq 2(d(\lambda)-\ep)\phi(H)x\phi(H)-
4\phi(H)BxB\phi(H)+K,
\end{equation*}
$K\in\PsiSc^{-\infty,2}(\Xb,\calC)$.
\end{thm}

\begin{rem}\label{rem:Mourre}
Again, the presence of bundles would make no difference, as is apparent
from the proof given below and in the appendix.

This theorem could be proved directly, without the use of the global
positive commutator estimate, \eqref{eq:Mourre-0}, but for notational
(and reference) reasons, it is easier to proceed via \eqref{eq:Mourre-0}.
\end{rem}

\begin{proof}
First, note that changing $B$ by any term $B'=xB'_1\in x\DiffSc^1(\Xb;\calC)$
changes the left hand side by
\begin{equation}\label{eq:subpr-8}
\phi(H)i[B',H]\phi(H)=\phi(H)iB'(H-\lambda)\phi(H)
-\phi(H)i(H-\lambda)B'\phi(H).
\end{equation}
Now, $\|\phi(H)(H-\lambda)\|\leq\delta'$ if $\supp\phi\subset
(\lambda-\delta',\lambda+\delta')$. Thus, multiplying \eqref{eq:subpr-8}
from the left and right by $\tilde\phi(H)$, $\tilde\phi\in\Cinf_c(\Real;[0,1])$
supported in $(\lambda-\delta',\lambda+\delta')$,
$\phit(H)iB'(H-\lambda)\phit(H)$ has the form $\phit(H) x T\phit(H)$, and
$\|T\|\leq \delta'\|B'\|$.
Hence, after this multiplication, both terms on the right
hand side of \eqref{eq:subpr-8} can be absorbed into $2(d(\lambda)-\ep)
\phit(H)x\phit(H)$ at the cost of increasing $\ep>0$ (which was arbitrary
to start with) by an arbitrarily
small amount. Thus, for each $\ep>0$ ther existence of a $\delta>0$
such that the estimate of the theorem holds
only depends on the indicial operators of $B$, even though
$\phi(H)i[B',H]\phi(H)$ is the same order as our leading term
$2(d(\lambda)-\ep)\phi(H)x\phi(H)$ (the many-body
calculus is not commutative even to top order!): the key being that this
commutator is small at the `characteristic variety'.

With $x^2D_x=x(xD_x)$, $A=\half(\chi(x) xD_x+(\chi(x) xD_x)^*)$,
\begin{equation*}\begin{split}
i[x^2D_x,H]&=i[x,H](xD_x)+xi[xD_x,H]=i[x,H]x^{-1}(x^2D_x)
+xi[xD_x,H]\\
&=-(x^2D_x)x(x^2D_x)+x^{1/2}i[xD_x,H]x^{1/2}+K',
\ K'\in\PsiSc^{-\infty,2}(\Xb,\calC).
\end{split}\end{equation*}
Multiplying through by $\phi(H)$ from both the left and the right,
the standard Mourre estimate,
\begin{equation}\label{eq:Mourre-0}
\phi(H)i[A,H]\phi(H))\geq 2(d(\lambda)-\ep)\phi(H)^2+K',
\end{equation}
$K'\in\PsiSc^{-\infty,1}(\Xb,\calC)$, proves the Theorem, since
$x^{1/2}\phi(H)c\phi(H)x^{1/2}=\phi(H)cx\phi(H)+K''$,
$K''\in\PsiSc^{-\infty,2}(\Xb,\calC)$ as $[\phi(H),x^{1/2}]
\in\PsiSc^{-\infty,3/2}(\Xb,\calC)$.
As indicated above, we briefly recall the proof of the Mourre estimate in the
appendix, since technically our setting is not covered e.g.\ by the
proof of Froese and Herbst \cite{FroMourre}, even though their proof
goes through without any significant changes.
\end{proof}

An equivalent, and for us more useful, version of this theorem is the
following result.

\begin{cor}\label{cor:Mourre-est}
Suppose $\lambda\in\Real\setminus\Lambda$.
For any $\ep>0$, there exist $R\in x\PsiSc^{2,0}(\Xb,\calC)$,
$K\in x^2\PsiSc^{0,0}(\Xb,\calC)$, such that
\begin{equation}\label{eq:Mourre-ch}
i[B,H]\geq 2(d(\lambda)-\ep)x-
4BxB+(H-\lambda)R+R^*(H-\lambda)+ K.
\end{equation}
\end{cor}

\begin{rem}
This corollary essentially states that the commutator $i[B,H]$ is positive,
modulo $BxB$, on the `characteristic variety',
i.e.\ where $H-\lambda$ vanishes.
Since this is a non-commutative setting (even to leading order), the vanishing
on the characteristic variety has to be written by allowing error terms
$(H-\lambda)R+R^*(H-\lambda)$: $[R,H-\lambda]$ has the same order as $R$!

Also, the inequality \eqref{eq:Mourre-ch} is understood as a quadratic
form inequality on $\dCinf(\Xb)$.
\end{rem}

\begin{proof}
Let $\phi\in\Cinf_c(\Real)$ be identically $1$ near $\lambda$, supported
in $(\lambda-\delta,\lambda+\delta)$, $\delta$ as above. Then
\begin{equation*}\begin{split}
i[B,H]=&i\phi(H)[B,H]\phi(H)+i(\Id-\phi(H))[B,H]\phi(H)\\
&+i\phi(H)[B,H](\Id-\phi(H))
+i(\Id-\phi(H))[B,H](\Id-\phi(H)),
\end{split}\end{equation*}
with a similar expansion for $x$. Since $\Id-\phi(H)=(H-\lambda)\phit(H)$,
$\phit(t)=(t-\lambda)^{-1}(1-\phi(t))$, so $\phit\in S^{-1}(\Real)$
is a classical symbol, $\phit(H)\in\PsiSc^{-2,0}(\Xb,\calC)$. Thus,
\begin{equation*}
(\Id-\phi(H))[B,H]\phi(H)=(H-\lambda)R',\ R'\in\PsiSc^{-\infty,1}(\Xb,\calC),
\end{equation*}
etc., proving the corollary.
\end{proof}

\begin{rem}
This Corollary is in fact equivalent to the statement of the theorem. For
assuming \eqref{eq:Mourre-ch}, multiplying by $\phi_0(H)$ from the left and
right, $\phi_0\in\Cinf_c(\Real;[0,1])$ identically $1$ near $\lambda$,
replaces the term $(H-\lambda)R$ by $\phi_0(H)(H-\lambda)R\phi_0(H)$,
and $\phi_0(H)(H-\lambda)$ has small norm, $\leq\delta_0$,
if $\supp\phi_0\subset(\lambda-\delta_0,\lambda+\delta_0)$.
Multiplying from the left and right by $\phi(H)$ then gives
Theorem~\ref{thm:Mourre-est}, since these small terms can then be absorbed as
$2(d(\lambda)-\ep-C\delta_0)\phi(H)x\phi(H)$.
\end{rem}

\section{Exponential decay}
Using the preceeding lemma and the global positive commutator estimate,
Theorem~\ref{thm:Mourre-est}, we can now prove the exponential
decay of non-threshold eigenfunctions. For this part of the paper,
we could adapt the proof of
Froese and Herbst \cite{FroExp} in Euclidean potential scattering,
as was done in \cite{Vasy:Propagation-2} in the geometric potential
scattering setting. However, to unify the paper, we focus on the approach
that will play a crucial role in the proof of unique continuation at
infinity. Nonetheless, the Froese-Herbst commutator
will play a role when $\alpha=0$ (in the notation of Lemma~\ref{lemma:FH-1}),
where
conjugated Hamiltonian is close to being self-adjoint (in fact, it
is, if $F=0$), so we will use $xD_x$ for a commutator estimate
in place of $\im P$.

\begin{thm}\label{thm:exp-decay}
Let $\lambda\in\Real\setminus\Lambda$, and suppose that
$\psi\in L^2_{\scl}(\Xb)$ satisfies $H\psi=\lambda\psi$. Then $e^{\alpha/x}\psi
\in L^2_{\scl}(\Xb)$ for all $\alpha\in\Real$ such that
$[\lambda,\lambda+\alpha^2]\cap\Lambda=\emptyset$.
\end{thm}

\begin{proof}
The proof is by contradiction. First note that $\psi\in\dCinf(\Xb)$ by
a result of \cite{Vasy:Bound-States} which only makes use of positive
commutator estimates whose proof is unchanged in this greater
generality. Let
\begin{equation*}
\alpha_1=\sup\{\alpha\in[0,\infty):\ \exp(\alpha/x)\psi\in L^2_{\scl}(\Xb)\},
\end{equation*}
and suppose that $[\lambda,\lambda+\alpha_1^2]\cap\Lambda=\emptyset$.
If $\alpha_1=0$, then let $\alpha=0$,
otherwise suppose that $\alpha<\alpha_1$,
and $\alpha+\gamma>\alpha_1$. We show that for sufficiently small $\gamma$
(depending only on $\alpha_1$) $\exp((\alpha+\gamma)/x)\psi\in L^2_{\scl}(\Xb)$,
which contradicts our assumption on $\alpha_1$ if $\alpha$ is close enough
to $\alpha_1$. In what follows we assume
that $\gamma\in(0,1]$.

Below we use two positivity estimates, namely \eqref{eq:FH-3} and the
Mourre-type
estimate, Corollary~\ref{cor:Mourre-est},
at energy $\lambda+\alpha_1^2$, with $B=\chi(x)x^2D_x+(\chi(x)x^2D_x)^*$.
That is, since
$\lambda+\alpha_1^2\nin\Lambda$,
there exists $c_0>0$,
$R\in\PsiSc^{0,0}(\Xb,\calC)$, $K\in\PsiSc^{2,0}(\Xb,\calC)$, such that
for $\psit\in L^2_{\scl}(\Xb)$,
\begin{equation}\begin{split}\label{eq:pos-comm-88}
(\psit,&i[B,H]\psit)\\
&\geq c_0\|x^{1/2}\psit\|^2-
4\re(\psit,x(x^2D_x)^2\psit)\\
&\qquad+\re((H-\lambda-\alpha_1^2)\psit,xR\psit)
+\re(x\psit,Kx\psit).
\end{split}\end{equation}
We apply this below with $\psit=\psi_F$.

We first note that we certainly have for all $\beta\in\Real$,
$\exp(\alpha/x)x^\beta\psi\in L^2_{\scl}(\Xb)$, due to our choice of $\alpha$.
We apply the Lemma~\ref{lemma:FH-1} with
\begin{equation*}
F=F_\beta=\frac{\alpha}{x}+\beta\log(1+\frac{\gamma}{\beta x}),
\end{equation*}
and let $\psi_\beta=e^F\psi$. (Since $x$ is bounded on $\Xb$, we may consider
$F$ compactly supported in $[0,x_0]$, $x_0=\sup_{\Xb}x<1/2$, as arranged
for convenience in the preceeding section.) The reason for this choice is that
on the one hand
$F(x)\to (\alpha+\gamma)/x$ as $\beta\to\infty$, so in the limit we will
obtain an estimate on $e^{(\alpha+\gamma)/x}\psi$, and on the other hand
$F(x)\leq \frac{\alpha}{x}+\beta|\log x|$, so $e^{F_\beta}$ is bounded
by $x^\beta e^{\alpha/x}$, for all values of $\beta$, i.e.\ $e^{F_\beta}$
provides a `regularization' (in terms of growth) of $e^{(\alpha+\gamma)/x}$,
so that Lemma~\ref{lemma:FH-1} can be applied.

Note that $F=F_\beta\in S^1([0,1))$, and
$F_\beta$ is uniformly bounded in $S^1(
[0,1))$
for $\beta\in[1,\infty)$, $\alpha\in[0,\alpha_1)$ (or $\alpha=\alpha_1$
if $\alpha_1=0$), $\gamma\in[0,1]$. Indeed,
\begin{equation*}
0\leq -x^2\pa_x F=\alpha+\gamma(1+\frac{\gamma}{\beta x})^{-1}\leq\alpha+\gamma,
\end{equation*}
and in general $(x\pa_x)^m(1+\frac{\gamma}{\beta x})^{-1}
=(-r\pa_r)^m(1+r)^{-1}$, $r=\frac{\gamma}{\beta x}$, so the uniform boundedness
of $F$ follows from $(1+r)^{-1}$ being a symbol in the usual sense on
$[0,\infty)$. In particular, all symbol norms of $-x^2\pa_x F-\alpha$ are
$\calO(\gamma)$. Below, when $\alpha=0$, we will need to consider
$(-x^2\pa_x F)^{-1}(x^2\pa_x)^m (-x^2\pa_x F)$, $m\geq 0$. By Leibniz' rule,
this can be written as
$\sum_{j\leq m}c_j x^{m} (-x^2\pa_x F)^{-1}(x\pa_x)^j (-x^2\pa_x F)$.
In terms of $r$, $(-x^2\pa_x F)^{-1}(x\pa_x)^j (-x^2\pa_x F)$ takes the form
$(1+r)(-r\pa_r)^m(1+r)^{-1}$, hence it
is still bounded on $[0,\infty)$, so in fact
\begin{equation}\label{eq:factor-8}
x^{-m}(-x^2\pa_x F)^{-1}(x^2\pa_x)^m (-x^2\pa_x F),\quad m\geq 0,
\end{equation}
is uniformly bounded on $[0,\infty)$. In fact, \eqref{eq:factor-8} is uniformly
bounded in $S^0([0,1))$, since applying $x\pa_x$ to it gives rise
to additional factors such as
\begin{equation*}
(-x^2\pa_x F)^{-k}(x\pa_x)^k(-x^2\pa_x F),
\end{equation*}
which are also uniformly bounded on $[0,\infty)$ by the same argument.

We remark first that $P(F)\psi_F=0$, so by elliptic regularity,
\begin{equation*}
\|\psi_F\|_{x^p\Hsc^k(\Xb)}\leq b_{1,k,p}\|x^p\psi_F\|,
\end{equation*}
with $b_{1,k,p}$ independent of $F$ as long as $\alpha$ is bounded.
{\em In general,
below $b_j$ denote positive
constants that are independent of $\alpha,\beta,\gamma$
in these intervals, and $R_j$ denote operators which are uniformly
bounded in $\DiffScc^2(\Xb,\calC)$, or on occasion in $\PsiScc^{m,0}
(\Xb,\calC)$, for some $m$.}
(Note that by elliptic regularity, the differential order
never matters.)

The proof is slightly different in the cases $\alpha>0$ and $\alpha=0$
since in the latter case the usually dominating term, $-2\alpha x^2D_x$,
of $\im P$ vanishes.

So assume first that $\alpha>0$.
The key step in the proof of this theorem arises from considering,
with $P=P_\beta=H(F)-\lambda$,
\begin{equation*}
P^*P=(\re P)^2+(\im P)^2+i(\re P\im P-\im P\re P),
\end{equation*}
so
\begin{equation}\label{eq:comm-8}
0=(\psi_F,P^*P\psi_F)=\|\re P\psi_F\|^2+\|\im P\psi_F\|^2+
(\psi_F,i[\re P,\im P]\psi_F).
\end{equation}
The first two terms on the right hand side are non-negative, so the
key issue is the positivity of the commutator. Note that
\begin{equation}\begin{split}\label{eq:comm-9}
\re P=H-\alpha^2-\lambda+\gamma R_1+xR_2,\\
\im P=-2\alpha x^2D_x+\gamma R_3+xR_4.
\end{split}\end{equation}
Below we use $H-\alpha^2-\lambda$ for
a positive commutator estimate, local in the spectrum of $H$, in place
of $\re P$, to make the choice
of the spectral cutoff $\phi$ independent of $\beta$ and $\gamma$.
(Otherwise we would need a uniform analogue of Theorem~\ref{thm:Mourre-est}
for $\re P$.)
Thus, by \eqref{eq:FH-3},
\begin{equation*}\begin{split}
i[\re P,\im P]=&2\alpha i[x^2D_x,H]+x\gamma R_5\\
&+xR_6\re P+\re PxR_7
+xR_8\im P+\im PxR_9
+x^2 R_{10},
\end{split}\end{equation*}
Hence, from \eqref{eq:comm-8} and \eqref{eq:pos-comm-88},
\begin{equation}\begin{split}\label{eq:comm-24}
0\geq &\|\re P\psi_F\|^2+\|\im P\psi_F\|^2+2\alpha c_0\|x^{1/2}\psi_F\|^2
+\gamma(\psi_F,xR_{11}\psi_F)\\
&+(\psi_F,xR_{12}\re P\psi_F)+(\psi_F,\re PxR_{13}\psi_F)\\
&+(\psi_F,xR_{14}\im P\psi_F)+(\psi_F,\im PxR_{15}\psi_F)
+(\psi_F,x^2R_{16}\psi_F).
\end{split}\end{equation}
Moreover, terms such as $|(\psi_F,x^2R_{16}\psi_F)|$ can be estimated
by $b_2\|x\psi_F\|^2$, while $\gamma|(\psi_F,xR_{11}\psi_F)|$
may be estimated by $\gamma b_3\|x^{1/2}\psi_F\|^2$, and
\begin{equation*}\begin{split}
&|(\psi_F,xR_{12}\re P\psi_F)|\leq b_4\|x\psi_F\|\|\re P\psi_F\|
\leq b_4(\ep^{-1}\|x\psi_F\|^2+\ep\|\re\psi_F\|^2)\\
&|(\psi_F,xR_{14}\im P\psi_F)|\leq b_5\|x\psi_F\|\|\im P\psi_F\|
\leq b_5(\ep^{-1}\|x\psi_F\|^2+\ep\|\im\psi_F\|^2)
\end{split}\end{equation*}
with similar estimates for the other terms.
Putting this together, \eqref{eq:comm-24} yields
\begin{equation}\begin{split}\label{eq:comm-32}
0\geq &(1-b_6 \ep)\|\re P\psi_F\|^2+(1-b_7 \ep)\|\im P\psi_F\|^2\\
&+(2\alpha c_0
-\gamma b_8)\|x^{1/2}\psi_F\|^2-b_9(\ep)\|x\psi_F\|^2.
\end{split}\end{equation}
For $\delta>0$, in $x\geq\delta$, $x|\psi_F|=xe^{F}|\psi|\leq
b_{10}(\delta)|\psi|$, so
\begin{equation*}\begin{split}
\|x\psi_F\|^2&=\|x\psi_F\|^2_{x\leq\delta}+\|x\psi_F\|^2_{x\geq\delta}\\
&\leq \delta\|x^{1/2}\psi_F\|^2_{x\leq\delta}
+b_{10}(\delta)\|\psi\|^2_{x\geq\delta}\\
&\leq \delta\|x^{1/2}\psi_F\|^2
+b_{10}(\delta)\|\psi\|^2.
\end{split}\end{equation*}
Thus, \eqref{eq:comm-32} yields that
\begin{equation}\begin{split}\label{eq:comm-64}
0\geq &(1-b_6 \ep)\|\re P\psi_F\|^2+(1-b_7 \ep)\|\im P\psi_F\|^2\\
&+(2\alpha c_0
-\gamma b_8-b_9(\ep)\delta)\|x^{1/2}\psi_F\|^2-b_{10}(\delta)\|\psi\|^2.
\end{split}\end{equation}
Hence, choosing $\ep>0$ sufficiently small so that $b_6\ep<1$, $b_7\ep<1$,
then choosing $\gamma_0>0$ sufficiently small so that $b_{11}=2\alpha c_0
-\gamma_0 b_8>0$, we deduce that for $\gamma<\gamma_0$,
\begin{equation}\label{eq:comm-88}
b_{10}(\delta)\|\psi\|^2\geq (b_{11}-b_9\delta)\|x^{1/2}\psi_F\|^2.
\end{equation}
But, for $\delta\in(0,\frac{b_{11}}{b_9})$,
this shows that $\|x^{1/2}\psi_F\|^2$
is uniformly bounded as $\beta\to\infty$.
Noting that $F$ is an increasing function of $\beta$
and $\psi_F$ converges to $e^{(\alpha+\gamma)/x}\psi$ pointwise,
we deduce from the monotone convergence theorem that
\begin{equation*}
x^{1/2}e^{(\alpha+\gamma)/x}\psi\in L^2_\scl(\Xb),
\end{equation*}
so for
$\gamma'<\gamma$, $e^{(\alpha+\gamma')/x}\psi\in L^2_\scl(\Xb)$.

In case $\alpha=0$, \eqref{eq:comm-8} still yields that, as
$i[\re P,\im P]=\gamma x R_{17}$, with $R_{17}$ uniformly bounded, that
$\|\re P\psi_F\|\leq b_{12}\|x^{1/2}\psi_F\|$,
$\|\im P\psi_F\|\leq b_{12}\|x^{1/2}\psi_F\|$. In particular, the
former implies that
\begin{equation}\label{eq:deg-re-est}
\|(H-\lambda)\psi_F\|\leq \gamma b_{13}\|\psi_F\|+b_{14}\|x^{1/2}\psi_F\|,
\end{equation}
while the latter yields that
\begin{equation}\label{eq:deg-im-est}
\|(x^2\pa_x F)x^2 D_x\psi_F\|\leq b_{15}\|x^{1/2}\psi_F\|.
\end{equation}
However, instead of the degenerating commutator $[\re P,\im P]$, we
consider $P^*A-AP$, with $A$ as in Theorem~\ref{thm:Mourre-est},
which is the expression considered by Froese
and Herbst in \cite{FroMourre}. Since
$A$ is $xD_x$, modulo lower order terms, and $\im P$ is $2(x^2\pa_x F)
(x^2D_x)$,
modulo lower order terms, $A$ can be considered a rescaling of $\im P$,
in that the degenerating factor $x^2\pa_x F$ is removed.
Now,
\begin{equation*}\begin{split}
&i(AP-P^*A)=i[A,\re P]-(\im P A+A\im P)\\
&\qquad=i[A,H-\lambda]+ xR_{18}
-4Ax(x^2\pa_x F)A+R_{19}xA+AxR_{20}+\gamma R_{21}.
\end{split}\end{equation*}
Now, by \eqref{eq:factor-8} and Lemma~\ref{lemma:FH-1},
$R_{20}=(x^2\pa_x F)R'_{20}+xR''_{20}$ with $R'_{20}, R''_{20}$
also bounded, and
similarly for $R_{19}$. Thus,
\begin{equation*}\begin{split}
0&=(\psi_F,i(AP-P^*A)\psi_F)\\
&\geq (\psi_F,i[A,H]\psi_F)
+4\|x^{1/2}(-x^2\pa_x F)^{1/2}A\psi_F\|^2\\
&\qquad- b_{16}\|x^{1/2}\psi_F\|^2
-b_{17}\|\psi_F\|\|(x^2\pa_x F) x^2D_x \psi_F\|-b_{18}\gamma\|\psi_F\|^2.
\end{split}\end{equation*}
Using the Mourre estimate
\eqref{eq:Mourre-0}, with $\phi(H)$ dropped but $H-\lambda$ inserted, as
in Corollary~\ref{cor:Mourre-est},
we deduce that (with $c'_0>0$)
\begin{equation*}\begin{split}
0\geq c'_0\|\psi_F\|^2+4&\|x^{1/2}(-x^2\pa_x F)^{1/2}A\psi_F\|^2-b_{21}
\|(H-\lambda)\psi_F\|\|\psi_F\|
-b_{16}\|x^{1/2}\psi_F\|^2\\
&\qquad\qquad
-b_{17}\|\psi_F\|\|(x^2\pa_x F) x^2D_x \psi_F\|-b_{18}\gamma\|\psi_F\|^2.
\end{split}\end{equation*}
Using \eqref{eq:deg-re-est}-\eqref{eq:deg-im-est} we deduce, as above, that
\begin{equation*}\begin{split}
0&\geq c'_0\|\psi_F\|^2-\gamma b_{21}\|\psi_F\|^2-b_{22}\|\psi_F\|\|x^{1/2}
\psi_F\|-b_{23}\|x^{1/2}\psi_F\|^2\\
&\geq (c'_0-\gamma b_{21}-\ep_1 b_{22})\|\psi_F\|^2-(b_{22}\ep_1^{-1}+b_{23})
\|x^{1/2}\psi_F\|^2\\
&\geq(c'_0-\gamma b_{21}-\ep_1 b_{22}-(b_{22}\ep_1^{-1}+b_{23})\delta)
\|\psi_F\|^2-b_{22}\ep_1^{-1}
\|\psi_F\|^2.
\end{split}\end{equation*}
Again, we fix first $\ep_1>0$ so that $c'_0-\ep_1 b_{22}>0$, then $\gamma_0>0$
so that $c'_0-\gamma_0 b_{21}-\ep_1 b_{22}>0$, finally $\delta>0$
so that $c'_0-\gamma_0 b_{21}-\ep_1 b_{22}-(b_{22}\ep_1^{-1}+b_{23})\delta>0$.
Now letting $\beta\to \infty$ gives that $e^{\gamma/x}\psi\in L^2_{\scl}(\Xb)$
for $\gamma<\gamma_0$,
as above.
\end{proof}

Having proved the exponential decay of non-threshold eigenfunctions, we can
also prove that the thresholds are isolated from above inductively, using
an observation of Perry \cite{Perry:Exponential}. This relies on the
following uniform estimate.

\begin{prop}\label{prop:unif-est}
Suppose that $\lambda_0\nin\Lambda$, and let $I$ be a compact interval
with $\sup I<\lambda_0$. Then there
exists $C>0$ with the following property. If $H\psi=\lambda\psi$, $\lambda
\in I$ and if $e^{\alpha'/x}\psi\in L^2_\scl(\Xb)$ for some $\alpha'>
\sqrt{\lambda_0-\lambda}$ then $\|x^{1/2}e^{\alpha/x}\psi\|\leq C\|\psi\|$
for $\alpha=\sqrt{\lambda_0-\lambda}$.
\end{prop}

\begin{proof}
The proof is very close to that of the preceeding theorem. First, we may
use $F=\alpha/x$ directly, i.e.\ take $\gamma=0$. Again, all constants are
uniform in $\alpha$ and $\psi$, provided that $\alpha$ is bounded.
Thus, \eqref{eq:comm-88} yields that
\begin{equation}\label{eq:comm-88p}
b_{10}(\delta)\|\psi\|^2\geq (b_{11}-b_9\delta)\|x^{1/2}\psi_F\|^2.
\end{equation}
So taking $\delta\in (0, \frac{b_{11}}{b_9})$ shows that
\begin{equation*}
\|x^{1/2}e^{\alpha/x}\psi\|^2\leq b_{12}'\|\psi\|^2,
\end{equation*}
which proves the proposition.
\end{proof}

We introduce the following terminology. If $H$ is a many-body Hamiltonian,
we say that the thresholds $\lambda\in\Lambda$ of
$H$ are isolated from above if
\begin{equation*}
\lambda\in \Lambda\Rightarrow
\exists\lambda'>\lambda\Mst(\lambda,\lambda')\cap
(\Lambda\cup\pspec(H))=\emptyset.
\end{equation*}

\begin{thm}
Let $\Lambda$ be the set of thresholds of $H$, and suppose that $\lambda\in
\Lambda$. Then $\lambda$ is isolated from above in $\Lambda\cup\pspec(H)$,
i.e.\ there exists $\lambda'>\lambda$ such that $(\lambda,\lambda')
\cap(\Lambda\cup\pspec(H))=\emptyset$.
\end{thm}

\begin{proof}
Note that the statement of the theorem is certainly true for $H_0$,
since $\Lambda_0\cup\pspec(H_0)=\{0\}$.
We prove inductively that if in all proper subsystems $H_b$ of $H_a$, the
thresholds are isolated from above, then the same holds for $H_a$.
Assuming the inductive hypothesis, and recalling
that
\begin{equation*}
\Lambda_a=\cup_{X^b\subsetneq X^a}
(\Lambda_b\cup\pspec(H^b)),
\end{equation*}
with the union being finite, we deduce that for any $\lambda\in
\Lambda_a$ there exists $\lambda'>\lambda$ such that $(\lambda,\lambda')
\cap\Lambda_a=\emptyset$. So we only need to show that $\pspec(H)
\cap(\lambda,\frac{\lambda+\lambda'}{2})$ is finite.

Suppose otherwise, and let $\psi^{(j)}$ be an orthonormal sequence of
eigenfunctions with eigenvalue $\lambda_j\in
(\lambda,\frac{\lambda+\lambda'}{2})$, $j\geq 1$.
Since by the Mourre
estimate eigenvalues may only accumulate at thresholds, $\lim_{j\to\infty}
\lambda_j=\lambda$. In particular, dropping $\psi_j$ for a finite number of
$j$, we may assume that $\lambda_j<\frac{7\lambda+\lambda'}{8}$ for all $j$.
Let $\lambda_0=\frac{3\lambda+\lambda'}{4}$. By Theorem~\ref{thm:exp-decay},
for all $j$,
$\psi_j\in e^{-\gamma/x}L^2_\scl(\Xb)$ for
$\gamma<\sqrt{\frac{\lambda'-\lambda}{2}}$. Note that this holds
in particular for some $\gamma>\sqrt{\lambda_0-\lambda}>
\sqrt{\lambda_0-\lambda_j}$ for all $j$.
Hence, by Proposition~\ref{prop:unif-est}, with
$\alpha_j=\sqrt{\lambda_0-\lambda_j}>
\frac{1}{4}\sqrt{\lambda'-\lambda}=\alpha_0$, there exists $C>0$
such that $\|x^{1/2}e^{\alpha_j/x}\psi_j\|_{\Hsc^2(\Xb)}\leq C$ for all $j$,
hence
$\|e^{\alpha_0/x}\psi_j\|_{\Hsc^2(\Xb)}\leq C'$ for all $j$. But
the inclusion $e^{-\alpha_0/x}\Hsc^2(\Xb)
\hookrightarrow L^2_{\scl}(\Xb)$ is compact, so $\psi_j$ has
a subsequence that converges in $L^2_{\scl}(\Xb)$,
which contradicts the orthogonality of the $\psi_j$. This completes
the inductive step, proving the theorem.
\end{proof}

\section{Absence of positive eigenvalues -- high energy estimates}
We next prove that faster than exponential
decay of an eigenfunction of $H$ implies that it vanishes. This was
also the approach taken by Froese and Herbst. However, we use a different,
more robust, approach to deal with our much larger error terms.
The proof is based on conjugation by $\exp(\alpha/x)$ and letting
$\alpha\to+\infty$. Correspondingly, we require positive commutator
estimates at high energies. In such a setting first order terms are
irrelevant, i.e.\ $V$ does not play a significant role below. On the
other hand, $\Delta_g-\Delta_{g_0}$ is not negligible in any sense.
However, we show that if $g$ and $g_0$ are close in a $\Cinf$ sense
(keeping in mind that we are assuming that $g$ has a special structure),
then the corresponding unique continuation theorem is still true.
Our argument also shows the very close connection with H\"ormander's
solvability condition. Indeed, we work semiclassically (writing
$h=\alpha^{-1}$), and the key fact we use is that the commutator
of the real and imaginary parts of the conjugated Hamiltonian has
the correct sign on its `non-commutative characteristic variety'.

\begin{thm}\label{thm:unique}
Let $\lambda\in\Real$ and let $d$ denote a metric giving the usual topology on
$\Cinf$ sections of $\ScT[\Xb;\calC]\otimes\ScT[\Xb;\calC]$.
There exists $\ep>0$ such that if
$d(g,g_0)<\ep$ and $H\psi
=\lambda\psi$, $\exp(\alpha/x)\psi\in L^2_{\scl}(\Xb)$ for all $\alpha$,
then $\psi=0$.
\end{thm}

\begin{proof}
Let $F=F_\alpha=\phi(x)\frac{\alpha}{x}$ where $\phi\in\Cinf_c(\Real)$ is
supported near $0$, identically $1$ in a smaller neighborhood of $0$,
and let $\psi_F=e^F\psi$. Then with $h=\alpha^{-1}$, $H_h=h^2H(F)$ and
$P_h=H_h-h^2\lambda$
are elliptic semiclassical differential operators, elliptic in the usual
sense (differentiability), and
\begin{equation*}
P_h\psi_h=0,\ \psi_h=\psi_F,
\end{equation*}
so by elliptic regularity,
\begin{equation}\label{eq:h-ell-reg}
\|\psi_h\|_{x^pH^k_{\scl,h}(\Xb)}\leq C_1\|\psi_h\|_{x^pL^2_{\scl}(\Xb)},
\end{equation}
$C_1$ independent of $h\in(0,1]$ (but depends on $k$ and $p$). In general,
below the $C_j$ denote constants independent of $h\in(0,1]$ (and $\delta>0$).

The key step in the proof of this theorem arises from considering
\begin{equation*}
P_h^*P_h=(\re P_h)^2+(\im P_h)^2+i(\re P_h\im P_h-\im P_h\re P_h)
\end{equation*}
so
\begin{equation}\label{eq:h-comm-8}
0=(\psi_h,P_h^*P_h\psi_h)=\|\re P_h\psi_h\|^2+\|\im P_h\psi_h\|^2+
(\psi_h,i[\re P_h,\im P_h]\psi_h).
\end{equation}
The first two terms on the right hand side are non-negative, so the
key issue is the positivity of the commutator. More precisely, we
need that there exist operators $R_j$ bounded in
$\Diff_{\Scl,h}^{2,0}(\Xb,\calC)$ such that
\begin{equation}\begin{split}\label{eq:h-comm-16}
(\psi_h&,i[\re P_h,\im P_h]\psi_h)\\
&\geq (\psi_h,
(xh+\re P_h xhR_1+xhR_2\re P_h\\
&\qquad\qquad+\im P_h xhR_3+xhR_4\im P_h+xh^2R_5+x^2hR_6)\psi_h).
\end{split}\end{equation}
The important
point is that replacing both $\re P_h$ and $\im P_h$ by zero, the
commutator is estimated from below by a positive multiple of $xh$,
plus terms $O(xh^2)$ and $O(x^2h)$.

We first prove \eqref{eq:h-comm-16}, and then show how to use it to
prove the theorem. First, modulo terms that will give contributions that
are in the error terms, $\re P_h$ may be replaced by $h^2\Delta_g-1$,
while $\im P_h$ may be replaced by $-2h(x^2D_x)$. Now, by a principal
symbol calculation (which also gives the `trivial case' of
Theorem~\ref{thm:Mourre-est}),
\begin{equation*}
i[h^2\Delta_{g_0}-1,-2h(x^2D_x)]=xh(4h^2\Delta_{g_0}-4h^2(x^2 D_x)^2+R_7),
\ R_7\in x\Diff^2_{\scl,h}(\Xb).
\end{equation*}
The key point here is the microlocal positivity of the commutator where
$h^2\Delta_{g_0}-1$ and $-2h (x^2D_x)$ both vanish.
Now, taking the commutator with $h x^2D_x$ is continuous from
$\Diff_{\Scl,h}^{2,0}(\Xb,\calC)$ to $xh\Diff_{\Scl,h}^{2,0}(\Xb,\calC)$,
so
\begin{equation*}
i[h^2\Delta_{g}-1,-2h(x^2D_x)]=xh(4h^2\Delta_{g}+R_8-4h^2(x^2 D_x)^2+R_7),
\end{equation*}
$R_7\in x\Diff^2_{\scl,h}(\Xb)$, $R_8\in \Diff^2_{\scl,h}(\Xb)$,
and
\begin{equation*}
\|R_8\|_{\bop(x^{1/2}H^2_{\scl,h}(\Xb),x^{1/2}L^2_{\scl}(\Xb))}\leq \rho(d(g,g_0)),
\end{equation*}
with $\rho$ continuous, $\rho(0)=0$.
Since
\begin{equation*}
\|x^{1/2}R_8\psi_h\|\leq
\|R_8\|_{\bop(x^{1/2}H^2_{\scl,h}(\Xb),x^{1/2}L^2_{\scl}(\Xb))}
\|\psi_h\|_{x^{1/2}H^2_{\scl,h}(\Xb)},
\end{equation*}
we deduce from \eqref{eq:h-ell-reg} that
\begin{equation*}
|(\psi_h,hxR_8\psi_h)|\leq h\|x^{1/2}\psi_h\|\|x^{1/2}R_8\psi_h\|
\leq C_1 h \rho(d(g,g_0))\|x^{1/2}\psi_h\|^2
\end{equation*}
This proves \eqref{eq:h-comm-16} if $C_1\rho(d(g,g_0))<3$, hence if $g$ is
close to $g_0$.

We now show how to use \eqref{eq:h-comm-16} to show unique
continuation at infinity. Let $x_0=\sup_{\Xb} x$.
We first remark that
\begin{equation*}\begin{split}
&|(\psi_h, xh R_2\re P_h\psi_h)|\leq C_2 h\|x\psi_h\|\|\re P_h\psi_h\|
\leq C_2 h\|x\psi_h\|^2+C_2 h\|\re P_h\psi_h\|^2,\\
&|(\psi_h, xh R_4\im P_h\psi_h)|\leq C_3 h\|x\psi_h\|\|\im P_h\psi_h\|
\leq C_3 h\|x\psi_h\|^2+C_3 h\|\im P_h\psi_h\|^2,
\end{split}\end{equation*}
with similar expressions for the $R_1$ and $R_3$ terms in
\eqref{eq:h-comm-16}. Next,
\begin{equation*}\begin{split}
&|(\psi_h,xh^2R_5\psi_h)|\leq C_4 h^2\|x^{1/2}\psi_h\|^2\\
&|(\psi_h,x^2h R_6\psi_h)|\leq C_5 h\|x\psi_h\|^2.
\end{split}\end{equation*}
For $\delta>0$, in $x\geq\delta$, $|\psi_h|=e^{1/{xh}}|\psi|\leq
e^{1/(\delta h)}|\psi|$, so
\begin{equation*}\begin{split}
\|x\psi_h\|^2&=\|x\psi_h\|^2_{x\leq\delta}+\|x\psi_h\|^2_{x\geq\delta}\\
&\leq \delta\|x^{1/2}\psi_h\|^2_{x\leq\delta}
+x_0^2 e^{2/(\delta h)}\|\psi\|^2_{x\geq\delta}\\
&\leq \delta\|x^{1/2}\psi_h\|^2
+x_0^2 e^{2/(\delta h)}\|\psi\|^2.
\end{split}\end{equation*}
Thus,
\begin{equation*}
\|(\psi_h,x^2h R_6\psi_h)|\leq C_5 h\delta\|x^{1/2}\psi_h\|^2
+C_5 x_0^2 he^{2/(\delta h)}\|\psi\|^2.
\end{equation*}
Hence, we deduce from \eqref{eq:h-comm-8}-\eqref{eq:h-comm-16} that
\begin{equation*}\begin{split}
0\geq (1-C_6 h)\|\re P_h\psi_h\|^2+(1-C_7 h)\|\im P_h\psi_h\|^2
&+h(1-C_8 h-C_9\delta)\|x^{1/2}\psi_h\|^2\\
&-C_{10}he^{2/(\delta h)}\|\psi\|^2.
\end{split}\end{equation*}
Hence, there exists $h_0>0$ such that for $h\in(0,h_0)$,
\begin{equation}\label{eq:h-comm-64}
C_{10} h e^{2/(\delta h)}\|\psi\|^2
\geq h(\frac{1}{2}-C_9\delta)\|x^{1/2}\psi_h\|^2
\end{equation}

Now suppose that $\delta\in (0,\min(\frac{1}{4C_9},\frac{1}{h_0}))$ and
$\supp\psi\cap\{
x\leq\frac{\delta}{4}\}$ is non-empty.
Since $xe^{2/xh}=h^{-1}f(xh)$ where $f(t)=te^{2/t}$, and $f$ is decreasing
on $(0,2)$ (its minimum on $(0,\infty)$ is assumed at $2$), we deduce
that for $x\leq\delta/2$, $x e^{2/xh}\geq \frac{\delta}{2h} e^{4/(\delta h)}$,
so
\begin{equation*}
\|x^{1/2}\psi_h\|^2\geq C_{11} \delta h^{-1} e^{4/(\delta h)},\ C_{11}>0.
\end{equation*}
Thus, we conclude from \eqref{eq:h-comm-64} that
\begin{equation*}
C_{10}\|\psi\|^2\geq (\frac{1}{2}-C_9\delta)C_{11} \delta h^{-1}
e^{2/(\delta h)}.
\end{equation*}
But letting $h\to 0$,
the right hand side goes to $+\infty$, providing a contradiction.

Thus, $\psi$ vanishes for $x\leq \delta/4$, hence vanishes
identically on $\Xb$ by the usual Carleman-type unique continuation theorem
\cite[Theorem~17.2.1]{Hor}.
\end{proof}

\appendix

\section{Proof of the Mourre estimate, \eqref{eq:Mourre-0}}
In this section we recall briefly how the Mourre estimate,
\eqref{eq:Mourre-0}, is proved, relying on a now standard
iterative argument for the indicial
operators that originated in this form in \cite{FroMourre}; see also
\cite{Derezinski-Gerard:Scattering}.
Namely, to prove \eqref{eq:Mourre-0}, one only
needs to show that for all $b$, the corresponding indicial operators
satisfy the corresponding inequality, i.e.\ that
\begin{equation}\label{eq:Mourre-1}
\phi(\hat H_b)i\widehat{[A,H]}_b\phi(\hat H_b)
\geq 2(d(\lambda)-\ep)\phi(\hat H_b)^2.
\end{equation}
(This means that the operators on the two sides, which are families of
operators
on $X^b$, depending on $(y_b,\xi_b)\in\sct_{C_b}\Xb$, satisfy the inequality
for all $(y_b,\xi_b)\in\sct_{C_b}\Xb$.)
It is convenient to assume that $\phi$ is identically $1$ near
$\lambda$; if \eqref{eq:Mourre-1} holds for such $\phi$, it holds
for any $\phi_0$ with slightly smaller support, as follows by multiplication
by $\phi_0(\hat H_b)$ from the left and right.

Note that for $b=0$ the estimate
certainly holds: it comes from the Poisson bracket formula
in the scattering calculus, or from a direct computation yielding
$i\widehat{[A,H]}_0=2\Delta_{g_0}$. Hence, if the the localizing factor
$\phi(\hat H_0)=\phi(|\xi|^2)$ is supported in
$(\lambda-\delta,\lambda+\delta)$ and $\lambda>0$,
then \eqref{eq:Mourre-1} holds even
with $d(\lambda)-\ep$ replaced with $\lambda-\delta$. Note that
$\lambda\geq d(\lambda)$, if $\lambda>0$, since $0$ is a threshold of $H$.
On the other hand, if $\lambda<0$, both sides of \eqref{eq:Mourre-1}
vanish for $\phi$ supported near $\lambda$, so the inequality holds
trivially.

In general, we may assume inductively that at all clusters $c$ with
$C_c\subsetneq C_b$, i.e.\ $X^b\subsetneq X^c$, \eqref{eq:Mourre-1} has
been proved with $\phi$ replaced by a cutoff $\phit$ and $\ep$ replaced
by $\ep'$, i.e.\ we may assume that for all $\ep'>0$ there exists $\delta'>0$
such that for all $c$ with $C_c\subsetneq C_b$, and for all
$\phit\in\Cinf_c(\Real;[0,1])$ supported in $\lambda-\delta,\lambda+\delta)$,
\begin{equation}\label{eq:Mourre-1-0}
\phit(\hat H_c)i\widehat{[A,H]}_c\phit(\hat H_c)
\geq 2(d(\lambda)-\ep)\phit(\hat H_c)^2.
\end{equation}
But these are exactly the indicial operators of
$\phit(\hat H_b)i\widehat{[A,H]}_b\phit(\hat H_b)$, so, as
discussed in \cite[Proposition~8.2]{Vasy:Propagation-Many},
\eqref{eq:Mourre-1} implies that
\begin{equation}\label{eq:Mourre-8}
\phit(\hat H_b)i\widehat{[A,H]}_b\phit(\hat H_b)
\geq 2(d(\lambda)-\ep')\phit(\hat H_b)^2+K_b,
\ K_b\in\PsiSc^{-\infty,1}(X^b,\calC^b).
\end{equation}
Recall that this implication relies on a square root construction in
the many-body calculus, which is particularly simple in this case.

Now, we first multiply \eqref{eq:Mourre-8}
through by $\phi(H)$ from both the left and the
right. Recall that we use coordinates $(w_b,w^b)$ on $X_b\oplus X^b$ and
$(\xi_b,\xi^b)$ are the dual coordinates.
We remark that $\hat H_b=|\xi_b|^2+H^b$, so
if $\lambda-|\xi_b|^2$ is not an eigenvalue of $H^b$,
then as $\supp\phi\to\{\lambda\}$,
$\phi(H^b+|\xi_b|^2)\to 0$ strongly, so as $K_b$ is compact,
$\phi(H^b+|\xi_b|^2)K_b\to 0$ in norm; in particular it
can be made to have norm smaller than $\ep'-\ep>0$. After multiplication
from both sides by $\phi_1(\hat H_b)$, with $\phi_1$ having even smaller
support, \eqref{eq:Mourre-1}
follows (with $\phi_1$ in place of $\phi$),
with the size of $\supp\phi_1$ a priori depending on $\xi_b$.
However, $i\phi_1(\hat H_b)\widehat{[A,H]}_b\phi_1(\hat H_b)$
is continuous in $\xi_b$ with values in bounded operators on $L^2(X^b)$,
so if \eqref{eq:Mourre-1} holds at one value
of $\xi_b$, then it holds nearby. Moreover, for large $|\xi_b|$ both sides
vanish as $\hat H_b=H^b+|\xi_b|^2$, with $H^b$ bounded below, so
the estimate is in fact uniform
if we slightly increase $\ep>0$.

In general, if $E$
denotes the projection on the $L^2$ eigenspace of $H^b$ at $\lambda
-|\xi_b|^2$, the argument just sketched works if we can replace
$\phi(\hat H_b)$ by $(\Id-E)\phi(\hat H_b)$, in particular, it suffices to
show that
\begin{equation}\label{eq:Mourre-8p}
i\phit(\hat H_b)\widehat{[A,H]}_b\phit(\hat H_b)
\geq 2(d(\lambda)-\ep')\phit(\hat H_b)^2+(\Id-E)K'_b(\Id-E),
\end{equation}
$K'_b$ compact on $L^2(X^b)$.

To show \eqref{eq:Mourre-8p}, we follow an argument due to B.\ Simon
(as explained
in the paper \cite{FroMourre} of Froese and Herbst). The key point is
to replace $E$ by a finite rank orthogonal projection $F$, which
will later ensure that an error term is finite rank, hence compact. Thus,
by the compactness of $K_b$ (from \eqref{eq:Mourre-8}),
there is a
finite rank orthogonal projection $F$ with $\Range F\subset\Range E$
(so $F$ commutes with $\hat H_b$)
such that
\begin{equation}\label{eq:F-proj-def}
\|(\Id-E)K_b(\Id-E)-(\Id-F)K_b(\Id-F)\|<\ep'.
\end{equation}
Multiplying \eqref{eq:Mourre-8} through by $(\Id-F)\phi(\hat H_b)
=\phi(\hat H_b)-F$ from
left and right and using \eqref{eq:F-proj-def} gives
\begin{equation}\begin{split}\label{eq:Mourre-16}
&i(\phi(\hat H_b)-F)\widehat{[A,H]}_b(\phi(\hat H_b)-F)\\
&\qquad\geq 2(d(\lambda)-2\ep')(\phi(\hat H_b)^2-F)+\phi(\hat H_b)(\Id-E)K_b
(\Id-E)\phi(\hat H_b).
\end{split}\end{equation}
Moving the terms involving $F$ from the left hand side to the right hand
side yields that
\begin{equation}\begin{split}\label{eq:Mourre-16p}
i\phi(\hat H_b)\widehat{[A,H]}_b\phi(\hat H_b)
&\geq 2(d(\lambda)-2\ep')(\phi(\hat H_b)^2-F)\\
&\qquad+F\phi(\hat H_b)i\widehat{[A,H]}_b(\phi(\hat H_b)-F)\\
&\qquad +(\phi(\hat H_b)-F)i\widehat{[A,H]}_b\phi(\hat H_b)F
+FTF\\
&\qquad+\phi(\hat H_b)(\Id-E)K_b
(\Id-E)\phi(\hat H_b).
\end{split}\end{equation}
where $T=Ei\widehat{[A,H]}_bE$.
But
\begin{equation*}\begin{split}
&F\phi(\hat H_b)i\widehat{[A,H]}_b(\phi(\hat H_b)-F)=FC+
F\phi(\hat H_b)i\widehat{[A,H]}_b(E-F),\\
&C=F\phi(\hat H_b)i\widehat{[A,H]}_b\phi(\hat H_b)(\Id-E),
\end{split}\end{equation*}
and for any $\ep_1>0$,
\begin{equation*}
FC+C^*F\geq-\ep_1 F-\ep_1^{-1} C^*C,
\end{equation*}
so \eqref{eq:Mourre-16p} yields that
\begin{equation}\begin{split}\label{eq:Mourre-32}
i\phi(\hat H_b)\widehat{[A,H]}_b\phi(\hat H_b)
\geq &2(d(\lambda)-2\ep')\phi(\hat H_b)^2\\
&-2(d(\lambda)-2\ep'+\ep_1/2)F\\
&+FT(E-F)+(E-F)TF+FTF\\
&
+\phi(\hat H_b)(\Id-E)(K_b+K''_b)
(\Id-E)\phi(\hat H_b),
\end{split}\end{equation}
with $K''_b=-\ep_1^{-1}
\widehat{[A,H]}_b F\phi(\hat H_b)^2 \widehat{[A,H]}_b$,
which is compact due to the appearance of $F$.
Now, $A=xD_x+A'=w_b\cdot D_{w_b}+w^b\cdot D_{w^b}+A'$,
$A'\in x\DiffSc^1(\Xb,\calC)$, hence, similarly
to \eqref{eq:subpr-8},
\begin{equation*}\begin{split}
T=Ei\widehat{[A,H]}_bE=&Ei\widehat{[w_bD_{w_b},H]}_bE+
Ei\widehat{[w^bD_{w^b},H]}_bE\\
&+E\widehat{A'}_b(\hat H_b-\lambda)E
-E(\hat H_b-\lambda)\widehat{A'}_b E=Ei\widehat{[w_b D_{w_b},H]}_bE.
\end{split}\end{equation*}
since by the virial theorem $iE[w^bD_{w^b},H^b]E=0$ (see the remark below).
Thus,
\begin{equation*}
T=iE\widehat{[w_b D_{w_b},\Delta_{X_b}]}E\geq2|\xi_b|^2 E=2(\lambda-\lambda')E\geq
2d(\lambda)E;
\end{equation*}
here $\lambda'=\lambda-|\xi_b|^2$ is the eigenvalue of $H^b$ to which
$E$ projects, and $\lambda-\lambda'\geq d(\lambda)$ since $\lambda'$ is a
threshold of $H$ by definition.
Note that $d(\lambda)$ enters the estimate at this point (i.e.\ this
is the constant we need to use, rather than
$\lambda$, which is the corresponding constant in the free region).
Thus,
\begin{equation*}
-2(d(\lambda)-2\ep'+\ep_1/2)F
+FT(E-F)+(E-F)TF+FTF\geq (4\ep'-\ep_1)F\geq 0
\end{equation*}
if we choose $\ep_1<4\ep'$.
Thus,
\begin{equation}\begin{split}\label{eq:Mourre-32p}
i\phi(\hat H_b)\widehat{[A,H]}_b\phi(\hat H_b)
&\geq 2(d(\lambda)-2\ep')\phi(\hat H_b)^2\\
&\qquad+\phi(\hat H_b)(\Id-E)(K_b+K''_b)
(\Id-E)\phi(\hat H_b),
\end{split}\end{equation}
which proves \eqref{eq:Mourre-8p}.
Hence the proof of \eqref{eq:Mourre-0} is complete.

\begin{rem}
The statement $iE[w^bD_{w^b},H^b]E=0$ is formally a consequence of
\begin{equation*}
iE[w^bD_{w^b},H^b]E=iEw^b D_{w^b}(H^b-\lambda-|\xi_b|^2)E
-iE(H^b-\lambda-|\xi_b|^2)w^b D_{w^b}E=0,
\end{equation*}
but this requires justification since $w^b D_{w^b}$ is not bounded
on $L^2(X^b)$. In fact, by elliptic regularity
(namely using $E=E\phi(\hat H_b)$), the only issue is
the lack of decay of $w^b D_{w^b}$ at infinity, but the computation is
justified by replacing $w^b D_{w^b}$ by $\chi_0(|w^b|/C) w^bD_{w^b}$,
$\chi_0\in\Cinf_c(\Real)$ identicaly $1$ near $0$, and observing that
\begin{equation*}
\phi(\hat H_b)[\chi_0(w^b/C),H]w^bD_{w^b}
\end{equation*}
is uniformly bounded and
goes to $0$ strongly as $C\to\infty$.
\end{rem}

\bibliographystyle{plain}
\bibliography{sm}
\end{document}